\documentclass[11pt]{article}
\usepackage{amssymb}
\usepackage{amsmath,amsfonts,amsthm,amscd,latexsym}
\def\uuar{\mathord{\mbox{\makebox[0pt][l]{\raisebox{.4ex}
{$\uparrow$}}$\uparrow$}}}

\def\dda{\mathord{\mbox{\makebox[0pt][l]{\raisebox{-.4ex}
{$\downarrow$}}$\downarrow$}}}

\def\da{\mathord{\downarrow}}

\def\ur{(U,R)}

\newtheorem{tm}{Theorem}[section]
\newtheorem{pn}[tm]{Proposition}
\newtheorem{lm}[tm]{Lemma}

\theoremstyle{definition}
\newtheorem{dn}[tm]{Definition}
\newtheorem{rk}[tm]{Remark}
\newtheorem{ex}[tm]{Example}

\usepackage[a4paper]{geometry}
\allowdisplaybreaks[4]
\begin{document}
\title{Representations of FS-domains and BF-domains via FS-approximation
Spaces\thanks{Supported by the Natural Science Foundation of China (11671008), the Natural Science Foundation of Jiangsu Province (BK20170483)}}
\author{{Guojun Wu}\\
{\small School of Mathematics and Statistics} \\ {\small Nanjing University of Information Science and Technology}\\
   {\small Nanjing 210044, P. R. China}\\
   {Luoshan Xu\thanks{Corresponding author. E-mail: luoshanxu@hotmail.com}}\\
{\small Department of Mathematics, Yangzhou
University, Yangzhou  225002, P.R. China}}

 \date{}
\maketitle

\begin{abstract}
In this paper, concepts of  (topological) FS-approximation
spaces are introduced.
Representations of FS-domains and BF-domains via (topological) FS-approximation
spaces are considered. It is proved that the collection of CF-closed sets in an FS-approximation
space (resp., a topological FS-approximation
space) endowed with the set-inclusion order is an FS-domain (resp., a BF-domain) and that every FS-domain (resp., BF-domain) is order  isomorphic to the  collection of CF-closed sets of some FS-approximation
space (resp., topological FS-approximation
space) endowed with the set-inclusion order. The concept of topological BF-approximation spaces is introduced and a skillful method without using CF-approximable relations to represent BF-domains is given. It is also proved  that the category of FS-domains (resp., BF-domains) with Scott continuous maps as morphisms is equivalent to that of FS-approximation spaces (resp., topological FS-approximation spaces) with CF-approximable relations as morphisms.\\

\indent \indent \indent \indent\  topological BF-approximation
space; Categorical equivalence
\end{abstract}

\section{Introduction}
Domain Theory,  developed from continuous lattices introduced by Scott \cite{Scott-72} in the  1970s as a denotational model for functional languages, is one of the important research fields of theoretical computer science  \cite{Gierz}. It was.  Mutual transformations and infiltration of  the mathematical structures of orders, topologies and algebras  are the basic features of this theory. In recent years, there is a growing body of scholarly work towards synthesizing Domain Theory with various new mathematical fields such as  Formal Concept Analysis \cite{Ganter}, Rough Set Theory \cite{Pawlak}, and Mathematical Logic.
Such syntheses are clearly reflected in the  research of  representations for various  domains.
By representations of domains, we mean any general way by which one can characterize a domain using a suitable family of some mathematical structures ordered by the set-inclusion order.
 There are  many ways to represent domains
via such as abstract bases \cite{wang2,xu-mao-ab}, formal topologies  \cite{XU-mao-form}, information systems  \cite{He-Xu,Spreen-xu-inf,wang2,Wang-Li,Xu-Mao2}, formal contexts \cite{guo1,guo2,guo3,wang,Wang-Wang-Li},  closure spaces \cite{Guo-Li,Li-Wang-Yao,Wang-Li-Guo,Wu-Guo-Li,Yao1-Li,Yao2-Li}, various kinds of logics \cite{Wang-log2,Wang-log1}.

Rough Set Theory, created by Pawlak and developed by many other mathematicians and computer scientists, is fundamentally
important in artificial intelligence and information sciences. It has provided a more general framework to express
common sense reasoning and uncertainty reasoning, and received wide attention on the research areas in both of
real-life applications and the theory itself.
Rough Set Theory  is closely related to Order Theory and Topology.
In \cite{jar2}, J\"{a}rvinen provided
the lattice-theoretical background of rough sets and studied order properties of generalized rough sets.  Yang and Xu in \cite{Yang-Xu-Alg} investigated algebraic and order structures of various families of subsets of
generalized approximation spaces (GA-spaces, for short) with possibly infinite universes and arbitrary binary relations. Inspired by Rough Set Theory, Zou, Li and Ho in  \cite{Zhou-Li-Ho} characterized continuous posets and Scott closed sets by approximation operators. In \cite{Lei-Luo}, Lei and Luo studied representations of complete lattices and algebraic lattices
based on Rough Set Theory. All of these reveal that there are deep relationships between Rough Set Theory and Domain Theory.

 Representations of domains using abstract bases appear to be the most natural and simple ones, while  the study scope of abstract bases is little narrow and easy to miss out on something deeper. Noticing from Rough Set Theory that abstract bases $(B, \prec)$ are  special GA-spaces $(U, R)$ \cite{Yang-Xu-Alg},  Wu and Xu in \cite{Wu-Xu} generalized  an abstract basis to a CF-approximation
space $(U, R, \mathcal{F})$, and generalized round ideals of abstract bases to CF-closed sets. Since  the lower approximation operator $\underline{R}$ and the upper approximation operator $\overline{R}$ are mutually dual, the upper approximation operator $\overline{R}$ was mainly used to give representations of domains via CF-approximation
spaces in \cite{Wu-Xu}. It turns out that this approach of  representing  domains is more general than that of representing  domains by abstract bases.  The concept of CF-approximable relations using a categorical approach was introduced and that the category of CF-approximation spaces and CF-approximable relations is equivalent to that of continuous domains and Scott continuous maps was proved.  The work in \cite{Wu-Xu} makes more natural and closed links  between Domain Theory and Rough Set Theory.

FS-domains were introduced and studied by  Jung in \cite{Jung1,Jung2,Jung3}. It is proved that the category of FS-domains (resp., BF-domains) is a maximal Cartesian closed full subcategory of the category of continuous (resp., algebraic) domains and Scott continuous maps.
  It is well known that a Cartesian closed category is of great importance in Domain Theory  as it can be employed to  model the typed $\lambda$-calculus. Based on this fact and the work in \cite{Wu-Xu},  this paper further works on  representations of FS-domains and BF-domains with special CF-approximation spaces---FS-approximation spaces and topological FS-approximation spaces.

This paper is organized as follows: In Section 2, we recall some  basic notions in Domain Theory and Rough Set Theory, as well as some   results in \cite{Wu-Xu}. In Section 3, we introduce the concept of FS-approximation spaces and discuss  properties of FS-approximation spaces. In Section 4, we give representations of FS-domains and BF-domains via (topological) FS-approximation spaces. In Section 5, we  will prove  that the category of FS-domains (resp., BF-domains) with Scott continuous maps as morphisms is equivalent to the category of FS-approximation spaces (resp., topological FS-approximation spaces) with CF-approximable relations as morphisms.
\section{Preliminaries}

In this section,   we quickly recall some basic notions and results in Domain Theory and Rough Set Theory. For  notions  not explicitly defined herein, the reader may refer to  \cite{Davey-Pri, Gierz,Jean-2013}. For some basic notions and results in category theory, please  refer to \cite{Well}.

For a set $U$ and $X\subseteq U$,  we use $\mathcal {P} (U)$
to denote the power set of $U$, $\mathcal{P}_{fin}(U)$ to denote the family of all nonempty finite subsets of $U$  and  $X^c$ to denote the complement
of $X$ in $U$. The symbol $F\subseteq_{fin} X$ means $F$ is a finite subset of $X$.

Let ($L$, $\leqslant$) be a poset. A \emph{principal ideal} (resp., \emph{principal filter}) of $L$
is a set of the form ${\da x}=\{y\in L\mid y\leqslant x\}$ (resp., ${\uparrow\! x}=\{y\in L\mid x\leqslant y\}$). For $A\subseteq L$, we write $\da A=\{y\in L\mid\exists \ x\in A, \,y\leqslant x\}$ and $\uparrow\! A=\{y\in L\mid\exists \ x\in A,\, x\leqslant
y\}$. A subset $A$ is a \emph{lower set} (resp., an \emph{upper set}) if $A=\da A$ (resp., $A=\uparrow\!A$). We say that $z$ is a \emph{lower bound} (resp., an \emph{upper bound}) of $A$ if $A\subseteq \uparrow\!z$ (resp., $A\subseteq\downarrow\!z$).  The supremum of $A$ is the least upper bound of $A$, denoted by $\bigvee A$ or $\sup A$. The infimum of $A$ is the greatest lower bound of $A$, denoted by $\bigwedge A$ or $\inf A$. A nonempty subset $D$ of $L$ is \emph{directed} if every  finite subset of $D$ has an upper bound in $D$. A subset $B$ of a directed set $D$ is called a cofinal subset of $D$, if for all $d\in D$ there is $a\in B$ such that $d\leqslant a$. It is well know that a cofinal subset of a directed set itself is a directed set.

 A poset $L$ is called a \emph{directed complete partially ordered set} ({\sl dcpo}, for short) if every directed subset of $L$ has a supremum.

%


\begin{lm}\label{lm-D=cup-direct} Let $D$ be a directed subset of a poset $P$ and $D$ has a supremum $\sup D$. If $\cup_{i=1}^n A_i=D$, then there is some $A_j \ (1\leqslant j\leqslant n)$ which is a cofinal subset of $D$, and $\sup A_j=\sup D$.
\end{lm}
\noindent{\bf Proof.} Assume that for all  $i\in\{1,2,\cdots, n\}$, $A_i$ is  not a cofinal subset of $D$. Then for every $i\in\{1,2,\cdots, n\}$, select $d_{i}\in D$ such that no element in $A_i$ greater than $d_{i}$. Since $D$ is directed and $\{d_{i}\mid 1\leqslant i\leqslant n\}$ is finite, there is $d\in D=\cup_{i=1}^n A_i$ such that  $\{d_{i}\mid 1\leqslant i\leqslant n\}\subseteq{\downarrow d}$. Hence, there is $j\in\{1,2,\cdots, n\}$ such that $d\in A_j$ and $d_{j}\leqslant d$,  contradicting to the choice of $d_j$. Thus, there is some $A_j \ (1\leqslant j\leqslant n)$ which is a cofinal subset of $D$.  Clearly,
 $\sup A_j=\sup D$.  \hfill$\Box$\\

Recall that in a poset $P$, we say that $x$ \emph{way-below}
$y$, written $x\ll y$,  if for any directed set $D$  having
a supremum with $\sup D\geqslant y$, there is some $d\in D$ such that $x\leqslant d$. If $x\ll x$, then $x$ is called a compact element of $P$. The set $\{x\in P\mid x\ll x\}$ is denoted by $K(P)$. 
The set $\{y\in P\mid  x\ll
y\}$ will be denoted $\uuar x$ and $\{y\in P\mid  y\ll x\}$ denoted
$\dda x$.  A poset $P$ is said to be
\emph{continuous} (resp., {\em algebraic}) if for all $ x\in P$, $\dda x$ is directed (resp., ${\downarrow} x\cap K(P)$ is directed) and  $x=\bigvee \dda x$ (resp., $x=\bigvee ({\downarrow} x\cap K(P))$). If a dcpo $P$ is  continuous (resp., algebraic), then $P$ is called a {\em continuous domain} (resp., {\em an algebraic domain}). A subset $B$ of of a poset $P$ is called a basis of $P$ if for all $x\in P$, there is $B_x\subseteq B\cap\dda x$ such that
$B_x$ is directed and $\sup B_x=x$. It is well known that a poset $P$ is continuous iff it has a
 basis,  and that $P$ is algebraic  iff $K(P)$ is a
 basis.
\begin{lm} {\em (\cite{Gierz})}\label{lm-way-rel} Let $P$ be a poset. Then for all $x, y, u, z\in P$,

{\em (1)} $x\ll y\Rightarrow x\leqslant y$;

{\em (2)} $u\leqslant x\ll y\leqslant z\Rightarrow u\ll z$.
\end{lm}
\begin{lm}{\em (\cite{Gierz})}\label{pn-int}
If $P$ is a continuous poset, then the way-below relation
$\ll$ has
the  interpolation property:
\vskip3pt

  \centerline{$x\ll z\Rightarrow \exists y\in P$ such that $x\ll y\ll
   z$.}
\end{lm}

Let $L$ and $P$ be dcpos, and $f: L \longrightarrow P$ a map. If for any directed subset $D\subseteq L$, $f(\bigvee D)=\bigvee f(D)$, then $f$ is called a {\em Scott continuous map}.
We denote by $[L\longrightarrow P]$ the poset of all the Scott continuous maps from $L$ to $P$ in the pointwise order: $f, g\in [L\longrightarrow P], f\leqslant g\Leftrightarrow f(x) \leqslant g(x)$ for all $x\in L$.

\begin{dn} (\cite{Gierz})\label{dn-doms}
Let  $L$ be a dcpo.

(1) {\em An approximate identity} for  $L$ is defined to be a directed set $\mathcal{D}\subseteq [L\longrightarrow L]$ satisfying $\sup\mathcal{D}=id_{L}$, where $id_{L}$ is the identity on $L$.

(2) A Scott continuous map $\delta: L\longrightarrow L$ is said to be {\em finitely separating} if there exists  a finite set $M_{\delta}$ such that for each $x\in L$, there exists $m\in M_{\delta}$
satisfying  $\delta(x)\leqslant m\leqslant  x$.

(3) If there is an approximate identity for $L$ consisting of finitely separating maps, then $L$ is called  {\em an FS-domain}.

(4) If $L$ is an algebraic FS-domain, then $L$ is called  {\em a BF-domain}.
\end{dn}

\begin{lm}{\em (\cite[Lemma II-2.14]{Gierz})}\label{lm-fis-map} Let L be a dcpo.

{\em(1)} If $\mathcal{D}\subseteq [L\longrightarrow L]$ is an approximate identity for $L$, then $\mathcal{D}^{\prime}=\{\delta^{2}=\delta\circ\delta\mid \delta\in \mathcal{D}\}$ is also an approximate identity.

{\em(2)} If $\delta\in [L\longrightarrow L]$ is finitely separating, then $\delta(x)\ll x$ for all $x \in  L$. Thus an FS-domain is a continuous domain.
\end{lm}

 The following definitions and lemmas are basic concepts and results in  Rough Set Theory.
\begin{dn} (\cite{gl-liu,Yang-Xu-Alg}) A
\emph{generalized approximation space}
(\emph{GA-space}, for short) is a pair $(U, R)$ consisting of a non-empty set $U$ and
a  non-empty binary relation $R$ on $U$. For a GA-space $\ur$,
define $R_s, R_p: U\rightarrow\mathcal {P}(U)$  such that for all $x\in U$,
\vskip3pt

\centerline{$R_s(x)=\{y\in U \mid   xRy\},\indent \indent R_p(x)=\{y\in U \mid  yRx\}.$}

\end{dn}

\begin{dn} (\cite{gl-liu,Yang-Xu-Alg})\label{dn-app}
 Let $(U,R)$ be a GA-space. Define $\underline{R},\overline{R}:\mathcal {P} (U)\rightarrow
\mathcal {P}(U)$ such that for $A\subseteq U$,
\vskip3pt

\centerline {$\underline{R} (A)=\{x\in U\mid \ R_s(x)\subseteq A\}$,\qquad
$\overline{R}(A)=\{x\in U\mid \ R_s(x)\cap A\neq \emptyset\}.$}
\vskip5pt

\noindent Operators $\underline{R}$ and $\overline{R}$  are called {\em the lower approximation operator} and {\em the upper
approximation operator} in $(U,R)$, respectively.
\end{dn}

\begin{lm} {\em (\cite{gl-liu,Yang-Xu-Alg})}\label{lm-ula}
Let $(U,R)$ be a GA-space.

{\em(1)} $\underline{R} (A^c)=(\overline{R}(A))^c$,
$\overline{R}(A^c)=(\underline{R}(A))^c$, where $A^c$ is the
complement of $A\subseteq U$.

{\em(2)} $\underline{R}(U)=U$, $\overline{R}(\emptyset)=\emptyset$.

{\em(3)} Let $\{A_{i}\mid i\in I\}\subseteq \mathcal {P}(U)$. Then
$\underline{R}(\bigcap_{i\in I}A_{i})=\bigcap_{i\in
I}\underline{R}(A_{i}),\ \overline{R}(\bigcup_{i\in
I}A_{i})=\bigcup_{i\in I}\overline{R}(A_{i}).$

{\em(4)} If $A\subseteq B\subseteq U$, then
$\underline{R}(A)\subseteq\underline{R}(B),\overline{R}(A)\subseteq\overline{R}(B)$.

{\em(5)} For all $x\in U, \ \overline{R}(\{x\})=R_p(x)$.
\end{lm}
 A binary relation $R\subseteq U\times U$ on a set $U$ is said {\em reflexive} if for all $x\in U$, one has $x R x$; and
is said \emph{transitive} if $xRy$ and $yRz$ implies $xRz$ for all
$x,y,z\in U$. A binary relation $R$ is called a {\em preorder} if it is  reflexive and transitive.
\begin{lm}\label{lm-ref-lm2-tr-up}{\em (\cite{Wu-Xu,yy-yao2})}
Let $(U,R)$ be a GA-space.

{\em(1)} $R$ is reflexive iff   for
all $X\subseteq U$,
 $X\subseteq \overline{R}(X)$.

{\em(2)} $R$ is transitive iff for
all $X\subseteq U$,
 $\overline{R}(\overline{R}(X))\subseteq \overline{R}(X)$.

{\em(3)} If $R$ is  transitive and $A, B\subseteq U$, then $\overline{R}(B)\subseteq\overline{R}(A)$
whenever $B\subseteq\overline{R}(A)$.

{\em(4)}  If $R$ is a  preorder, then the operator $\underline{R}$ is an interior operator of a topology on $U$.
\end{lm}

In \cite{Wu-Xu}, Wu and Xu introduced CF-approximation spaces and gave representations of continuous domains by them. We here recall some key notions and results  in \cite{Wu-Xu}.

\begin{dn}\label{dn-cf-ga}(\cite{Wu-Xu}) Let $(U, R)$ be a GA-space, $R$  a transitive relation and $\emptyset\ne\mathcal{F}\subseteq \mathcal{P}_{fin}(U)\cup\{\emptyset\}$. If for all $ F\in \mathcal{F}$,  whenever $K\subseteq_{fin} \overline{R}(F)$, there always exists $G\in \mathcal{F}$ such that $K\subseteq \overline{R}(G)$ and $G\subseteq \overline{R}(F)$, then $(U, R, \mathcal{F})$ is called {\em a generalized approximation space with consistent family of finite subsets}, or {\em a CF-approximation space}, for short.
\end{dn}
\begin{dn}\label{dn-cf-cl}  (\cite{Wu-Xu})  Let $(U, R, \mathcal{F})$ be a CF-approximation space, $E\subseteq U$. If for all $ K\subseteq_{fin}E$, there always exists $F\in \mathcal{F}$ such that $K\subseteq \overline{R}(F)\subseteq E$ and $F\subseteq E$, then $E$ is called {\em a CF-closed set} of $(U, R, \mathcal{F})$. The collection of all CF-closed sets of $(U, R, \mathcal{F})$ is denoted by $\mathfrak{C}(U, R, \mathcal{F})$.
\end{dn}

\begin{lm}{\em (\cite{Wu-Xu})}\label{pn2-cf-clo} For  a CF-approximation space $(U, R, \mathcal{F})$, the following statements hold:

$(1)$ For any $F\in \mathcal{F}$, $\overline{R}(F)\in \mathfrak{C}(U, R, \mathcal{F})$;

$(2)$ If $E\in \mathfrak{C}(U, R, \mathcal{F})$, $A\subseteq E$, then $\overline{R}(A)\subseteq E$;

$(3)$ If $\{E_{i}\}_{i\in I}\subseteq\mathfrak{C}(U, R, \mathcal{F})$ is a directed family, then $\bigcup_{i\in I}E_{i}\in \mathfrak{C}(U, R, \mathcal{F})$.
\end{lm}
\noindent
{\bf Proof.} (1) Let $F\in \mathcal{F}$. Then For all $ K\subseteq_{fin}\overline{R}(F)$, by Definition \ref{dn-cf-ga}, there is $G\in \mathcal{F}$ such that $K\subseteq \overline{R}(G)$ and $G\subseteq \overline{R}(F)$.  By Lemma \ref{lm-ref-lm2-tr-up}(3), $\overline{R}(G)\subseteq \overline{R}(F)$. By Definition \ref{dn-cf-cl}, we see that $\overline{R}(F)\in \mathfrak{C}(U, R, \mathcal{F})$.

(2) and (3)  follow  directly from  Lemma \ref{lm-ula}(4), \ref{lm-ref-lm2-tr-up}(3) and Definition \ref{dn-cf-cl}.\hfill$\Box$

\begin{lm}{\em (\cite{Wu-Xu})}\label{pn3-cf-cl}
For  a CF-approximation space $(U, R, \mathcal{F})$ and $E\subseteq U$, the following statements are equivalent:

$(1)$ $E\in \mathfrak{C}(U, R, \mathcal{F})$;

$(2)$ The family $\mathcal{A}=\{\overline{R}(F)\mid F\in \mathcal{F},  F\subseteq E\}$ is directed and $E=\bigcup\mathcal{A}$;

$(3)$ There exists a family $\{F_{i}\}_{i\in I}\subseteq \mathcal{F}$ such that $\{\overline{R}(F_{i})\}_{i\in I}$ is directed, and $E=\bigcup_{i\in I}\overline{R}(F_{i})$;

$(4)$  There always  exists $F\in \mathcal{F}$ such that $K\subseteq \overline{R}(F)\subseteq E$ whenever $ K\subseteq_{fin}E$.
\end{lm}
\noindent {\bf Proof.} If $E=\emptyset\in  \mathfrak{C}(U, R, \mathcal{F})$,  then by Definition \ref{dn-cf-cl}, $\emptyset\in \mathcal{F}$ and the lemma holds. If  $\emptyset\neq E\in \mathfrak{C}(U, R, \mathcal{F})$, then we have

$(1)\Rightarrow (2)$ By  Definition \ref{dn-cf-cl}, we know that $\mathcal{A}\ne\emptyset$. Let $X_{1}, X_{2}\in \mathcal{A}$, then there exist $F_{1}, F_{2}\in \mathcal{F}$ and $F_{1}, F_{2}\subseteq E$, such that $X_{1}=\overline{R}({F_{1}})$, $X_{2}=\overline{R}({F_{2}})$. By $F_{1}\cup F_{2}\subseteq_{fin} E$ and Definition \ref{dn-cf-cl} we know that there exists $F_{3}\in \mathcal{F}$, such that $F_{1}\cup F_{2}\subseteq\overline{R}(F_{3})$ and $F_{3} \subseteq E$.
 By  Lemma \ref{lm-ref-lm2-tr-up}(3), we know that $\overline{R}(F_{1})\subseteq\overline{R}(F_{3})$ and $\overline{R}(F_{2})\subseteq\overline{R}(F_{3})$. This shows that $\mathcal{A}$ is directed. Next we prove $E=\bigcup\mathcal{A}$. By Lemma \ref{pn2-cf-clo}(2) we know that $\bigcup\mathcal{A}\subseteq E$ holds. Conversely, if $x\in E$, then by Definition \ref{dn-cf-cl}, there is $F\in \mathcal{F}$ such that $x\in \overline{R}(F)\subseteq E$ and $F\subseteq E$. So $x\in \bigcup\mathcal{A}$. By the arbitrariness of $x\in E$ we know that $E\subseteq\bigcup\mathcal{A}$. Thus $E=\bigcup\mathcal{A}$.

$(2)\Rightarrow (3)$ Trivial.

$(3)\Rightarrow (4)$ It follows directly from  the finiteness of $ K$ and the directedness of  $\{\overline{R}(F_{i})\}_{i\in I}$.

$(4)\Rightarrow (1)$   If $ K\subseteq_{fin}E$, then there exists $F\in \mathcal{F}$ such that $K\subseteq \overline{R}(F)\subseteq E$.  By  Definition \ref{dn-cf-ga}, there exists $G\in \mathcal{F}$ such that $K\subseteq \overline{R}(G)$ and $G\subseteq \overline{R}(F)$. By Lemma \ref{lm-ref-lm2-tr-up}(3), we know that  $\overline{R}(G)\subseteq\overline{R}(F)\subseteq E$. Thus $K\subseteq \overline{R}(G)\subseteq E$. Noticing that $G\subseteq \overline{R}(F)\subseteq E$,  by Definition \ref{dn-cf-cl}, we have that $E\in \mathfrak{C}(U, R, \mathcal{F})$. \hfill$\Box$\\

\begin{ex} \label{tm-CDO-CFGA} Let $(L, \leqslant)$ be a continuous domain, $R_{L}$  the way-below relation ``$\ll$" of $(L, \leqslant)$; $\mathcal{F}_{L}=\{F\subseteq_{fin} L\mid F\  \mbox{has a top element } c_F\}$.

$(1)$ $(L, R_{L},  \mathcal{F}_{L})$ is a CF-approximation space.

$(2)$ $\mathfrak{C}(L, R_{L},  \mathcal{F}_{L})=\{\dda x\mid x\in L\}$.
\end{ex}
\noindent{\bf Proof.}  (1)   By Lemma \ref{lm-way-rel}, we know that $R_{L}=\ll$ is transitive. For any $F\in \mathcal{F}_{L}$, let $c_F$ be the top element in $F$. By Lemma \ref{lm-ula}(5), we have that $\overline{R}_{L}(F)=\dda c_{F}$.  For $K\subseteq_{fin}\overline{R}_{L}(F)=\dda c_{F}$, by that $L$ is a continuous domain, we know that $\dda c_{F}$ is directed. Then there exists $x\in \dda c_{F}$ such that $K\subseteq {\downarrow} x$. It follows from  $x\ll c_{F}$ and Lemma \ref{pn-int} that there is $y\in L$ such that $x\ll y\ll c_{F}$. Thus $K\subseteq \dda y$. Set $G=\{y\}\in \mathcal{F}_{L}$. By that $K\subseteq \overline{R}_{L}(G)=\dda y$ and $G\subseteq\overline{R}_{L}(F)=\dda c_{F}$, we have that $(L, R_{L},  \mathcal{F}_{L})$ is a CF-approximation space.

(2) By the proof of (1) and Lemma \ref{pn2-cf-clo}(1), we know that $\{\dda x\mid x\in L\}\subseteq \mathfrak{C}(L, R_{L},  \mathcal{F}_{L})$. Conversely, let $E\in \mathfrak{C}(L, R_{L},  \mathcal{F}_{L})$. Then by Lemma \ref{pn3-cf-cl} and the continuity of $L$, there is  a directed set $D\subseteq L$ such that $E=\bigcup\{\dda d\mid d\in D\}$.  Next we prove $E=\dda \bigvee D$. Obviously, $E\subseteq \dda \bigvee D$. Conversely, if $x\in \dda \bigvee D$, then by Lemma \ref{pn-int},  there is  $y\in L$ such that $x\ll y\ll\bigvee D$. So there is $d\in D$ such that $x\ll y\leqslant d$. Thus $x\in \dda d\subseteq E$, and  $E=\dda \bigvee D$. This shows that $\mathfrak{C}(L, R_{L},  \mathcal{F}_{L})\subseteq\{\dda x\mid x\in L\}$ and $\mathfrak{C}(L, R_{L},  \mathcal{F}_{L})=\{\dda x\mid x\in L\}$. \hfill$\Box$\\

 We call $(L, R_{L},  \mathcal{F}_{L})$  the induced CF-approximation space by continuous domain $(L, \leqslant)$.

\begin{lm}{\em (\cite{Wu-Xu})}\label{tm-cfga-way} Let $(U, R, \mathcal{F})$ be a CF-approximation space.

 $(1)$ If $E_{1}, E_{2}\in \mathfrak{C}(U, R, \mathcal{F})$,  then $ E_{1}\ll E_{2}$ if and only if there exists $F\in \mathcal{F}$ such that $ E_{1}\subseteq \overline{R}(F)$ and  $F\subseteq E_{2}$;

$(2)$ If $F\subseteq E$, then $\overline{R}(F)\ll E$;

$(3)$ $\overline{R}(F)\ll \overline{R}(F)$ if and only if there exists $G\in \mathcal{F}$, such that $G\subseteq\overline{R}(G)=\overline{R}(F)$.
\end{lm}
\noindent{\bf Proof.}
(1)  Let $E_{1}, E_{2}\in \mathfrak{C}(U, R, \mathcal{F})$. If $ E_{1}\ll E_{2}$, then it follows from $E_{2}\in \mathfrak{C}(U, R, \mathcal{F})$ and  Lemma \ref{pn3-cf-cl}(2)  that $E_{2}=\bigcup\{\overline{R}(F)\mid F\in \mathcal{F}, F\subseteq E_{2}\}$  and that $\{\overline{R}(F)\mid F\in \mathcal{F}, F\subseteq E_{2}\}$ is directed. By that $ E_{1}\ll E_{2}$,  there exists $F\in \mathcal{F}$ such that $F\subseteq E_{2}$ and $ E_{1}\subseteq \overline{R}(F)$. Conversely, if there exists $F\in \mathcal{F}$ such that $ E_{1}\subseteq \overline{R}(F)$ and  $F\subseteq E_{2}$, then for any directed family $\{C_{i}\}_{i\in I}\subseteq\mathfrak{C}(U, R, \mathcal{F})$ with $ E_{2}\subseteq \bigvee_{i\in I}C_{i}= \bigcup_{i\in I}C_{i}$, then by $F\subseteq E_{2}$ and the finiteness of $F$, we know that there exists $i_{0}\in I$ such that $F\subseteq C_{i_{0}}$. By Lemma \ref{pn2-cf-clo}(2)  and $E_{1}\subseteq \overline{R}(F)$, we know that $E_{1}\subseteq \overline{R}(F)\subseteq C_{i_{0}}$, showing that $E_{1}\ll E_{2}$.

(2) and (3)  follow  directly from  (1).\hfill$\Box$


\section{FS-approximation
spaces}

In this section, we introduce FS-approximation spaces and  topological FS-approximation spaces. They are   special types of CF-approximation spaces.  And then discuss their properties.
 We first recall the concept of CF-approximable relations  in \cite{Wu-Xu}.
 \begin{dn} (\cite{Wu-Xu}) \label{dn-CF-app}  Let $(U_{1}, R_{1}, \mathcal{F}_{1})$, $(U_{2}, R_{2}, \mathcal{F}_{2})$ be CF-approximation spaces, and $\mathrel{\Theta}\subseteq\mathcal{F}_{1}\times\mathcal{F}_{2}$ a binary relation.  If

(1)  for all $F\in \mathcal{F}_{1}$, there is $G\in \mathcal{F}_{2}$ such that $F\mathrel{\Theta} G$;

(2)  for all $F, F^{\prime}\in \mathcal{F}_{1}$, $G\in  \mathcal{F}_{2}$, if $F\subseteq \overline{R_{1}}(F^{\prime})$, $F\mathrel{\Theta} G$, then $F^{\prime}\mathrel{\Theta} G$;

(3)  for all $F\in \mathcal{F}_{1}$, $G, G^{\prime}\in  \mathcal{F}_{2}$, if $F\mathrel{\Theta} G$, $G^{\prime}\subseteq \overline{R_{2}}(G)$,  then $F\mathrel{\Theta} G^{\prime}$;

(4) for all $F\in \mathcal{F}_{1}$,  $G\in \mathcal{F}_{2}$, if $F\mathrel{\Theta} G$, then there are $F^{\prime}\in \mathcal{F}_{1}$,  $G^{\prime}\in \mathcal{F}_{2}$ such that  $F^{\prime}\subseteq \overline{R_{1}}(F)$, $G\subseteq \overline{R_{2}}(G^{\prime})$ and $F^{\prime}\mathrel{\Theta} G^{\prime}$; and

(5)  for all $F\in \mathcal{F}_{1}$, $G_{1}, G_{2}\in \mathcal{F}_{2}$, if $F\mathrel{\Theta} G_{1}$ and $F\mathrel{\Theta} G_{2}$, then there is  $G_{3}\in \mathcal{F}_{2}$ such that  $G_{1}\cup G_{2}\subseteq\overline{R_{2}}(G_{3})$ and $F\mathrel{\Theta} G_{3}$,

\noindent then $\mathrel{\Theta}$ is called {\em a CF-approximable relation} from $(U_{1}, R_{1}, \mathcal{F}_{1})$ to $(U_{2}, R_{2}, \mathcal{F}_{2})$.
\end{dn}
 \begin{dn}\label{dn-cf-app-rel-id}
Let $(U, R, \mathcal{F})$ be a CF-approximation space and
\vskip3pt

\centerline{$\operatorname{Id}_{(U, R, \mathcal{F})}=\{(F,  G)\mid F, G\in \mathcal{F}, G\subseteq\overline{R}(F)\} \subseteq\mathcal{F}\times\mathcal{F}$.}
\vskip5pt

\noindent
Then it is easy to check that $\operatorname{Id}_{(U, R, \mathcal{F})}$ is a CF-approximable relation from $(U, R, \mathcal{F})$ to itself. $\operatorname{Id}_{(U, R, \mathcal{F})}$ is called {\em the identity  CF-approximable relation} on $(U, R, \mathcal{F})$.
\end{dn}

\begin{pn} \label{pn1-CF-apprel} Let $\mathrel{\Theta}$ be a CF-approximable relation from  $(U_{1}, R_{1}, \mathcal{F}_{1})$ to  $(U_{2}, R_{2}, \mathcal{F}_{2})$. Then for all $ F\in \mathcal{F}_{1}, G\in \mathcal{F}_{2}$, the following statements are equivalent:

$(1)$ $F\mathrel{\Theta} G$;

$(2)$ There exists $F^{\prime}\in \mathcal{F}_{1}$ such that $F^{\prime}\subseteq \overline{R_{1}}(F)$ and $F^{\prime}\mathrel{\Theta} G$;

$(3)$ There exists $G^\prime\in \mathcal{F}_{2}$ such that $F\mathrel{\Theta} G^\prime$ and $G\subseteq\overline{R_{2}}(G^\prime)$;

$(4)$ There exist $F^{\prime}\in \mathcal{F}_{1}$ and  $G^{\prime}\in \mathcal{F}_{2}$ such that $F^{\prime}\subseteq \overline{R_{1}}(F)$, $G\subseteq \overline{R_{2}}(G^{\prime})$ and $F^{\prime}\mathrel{\Theta} G^{\prime}$.
\end{pn}
\noindent{\bf Proof. } Follows  directly from (2)-(4) in Definition \ref{dn-CF-app}.\hfill$\Box$
\begin{pn}\label{pn-ovRF-direct} Let $\mathrel{\Theta}$ be  a   CF-approximable relation from     $(U_{1}, R_{1}, \mathcal{F}_{1})$ to  $(U_{2}, R_{2}, \mathcal{F}_{2})$ and $E\in \mathfrak{C}(U_{1},  R_{1}, \mathcal{F}_{1})$. Then the family $\mathcal{D}=\{\overline{R}(G)\mid F\in \mathcal{F}_{1}, F\subseteq E, G\in \mathcal{F}_{2}\mbox{\ and\ } F\mathrel{\Theta}  G\}$  is directed.
 \end{pn}
\noindent{\bf Proof.}
 It is easy to see by Definition \ref{dn-cf-cl} and the condition (1) in Definition \ref{dn-CF-app} that $\mathcal{D}\ne \emptyset$. Let $X_{1}, X_2\in \mathcal{D}$. Then there are $F_{i}\in \mathcal{F}_{1}$ and $G_{i}\in \mathcal{F}_{2}$ such that $F_{i}\subseteq E$, $ F_{i}\mathrel{\Theta}  G_{i}$ and $X_{i}=\overline{R}(G_{i})$ $(i=1, 2)$. It follows from Definition \ref{dn-cf-cl} that there exits $F_{3}\in \mathcal{F}_{1}$ such that $F_{1}\cup F_{2}\subseteq \overline{R}(F_{3})$ and $F_{3}\subseteq E$.
It follows from Definition \ref{dn-CF-app}(2) that $F_{3}\mathrel{\Theta}  G_{1}$ and $F_{3}\mathrel{\Theta}  G_{2}$. By Definition \ref{dn-CF-app}(5), there exists $G_{3}\in \mathcal{F}_{2}$ such that  $G_{1}\cup G_{2}\subseteq\overline{R_{2}}(G_{3})$ and $F_{3}\mathrel{\Theta} G_{3}$. By Lemma \ref{lm-ref-lm2-tr-up} (3), we have that $X_{1}\cup X_{2}\subseteq \overline{R_{2}}(G_{3})$ and the family $\mathcal{D}=\{\overline{R}(G)\mid F\in \mathcal{F}_{1}, F\subseteq E, G\in \mathcal{F}_{2}\mbox{\ and\ } F\mathrel{\Theta}  G\}$ is directed. \hfill$\Box$

\begin{lm}{\em (\cite{Wu-Xu})}\label{tm-cfap-sco}
 Let $\mathrel{\Theta}$ be  a     CF-approximable relation from     $(U_{1}, R_{1}, \mathcal{F}_{1})$ to  $(U_{2}, R_{2}, \mathcal{F}_{2})$.  Define a map $f_{\mathrel{\Theta}}: \mathfrak{C}(U_{1}, R_{1}, \mathcal{F}_{1})\longrightarrow \mathfrak{C}(U_{2}, R_{2}, \mathcal{F}_{2})$ such that for all $ E\in \mathfrak{C}(U_{1}, R_{1}, \mathcal{F}_{1})$,
  $$f_{\mathrel{\Theta}}(E)=\bigcup\{\overline{R}(G)\mid F\in \mathcal{F}_{1}, F\subseteq E, G\in \mathcal{F}_{2}\mbox{\ and\ } F\mathrel{\Theta}  G\}.$$
Then  $f_{\mathrel{\Theta}}$ is a  Scott continuous map.
\end{lm}
\noindent {\bf Proof.}  It follows from Lemma \ref{pn2-cf-clo}, Proposition \ref{pn-ovRF-direct} and the definition of $f_{\mathrel{\Theta}}$ that $f_{\mathrel{\Theta}}$ is well-defined and order preserving. It is direct to show that for any directed family $\{E_{i}\}_{i\in I}\subseteq \mathfrak{C}(U_{1}, R_{1}, \mathcal{F}_{1})$, we have
$f_{\mathrel{\Theta}}(\bigcup_{i\in I}E_{i})=\bigcup_{i\in I}f_{\mathrel{\Theta}}(E_{i})$. By  Lemma \ref{pn2-cf-clo}(3), $f_{\mathrel{\Theta}}$ is a  Scott continuous map.\hfill$\Box$\\

Now we can give one of the main concepts in this paper.
\begin{dn}\label{dn-FS-ga} A CF-approximation space $(U, R, \mathcal{F})$ is called {\em an FS-approximation space} if there exists a directed family $\{\Theta_{i}\}_{i\in I}$ of  CF-approximable relations on $(U, R, \mathcal{F})$ satisfying the following conditions:

 (FS 1) $\bigcup_{i\in I}\Theta_{i}=\operatorname{Id}_{(U, R, \mathcal{F})}$;

(FS 2) For all $\Theta_{i}\ (i\in I)$, there is  $\mathcal{M}_{i}\subseteq_{fin}\mathcal{F}$ such that for all $F\in \mathcal{F}$,  exists  $M\in \mathcal{M}_{i}$ satisfying
$$\forall G\in \mathcal{F}, F\mathrel{\Theta_{i}}G \Rightarrow G\subseteq \overline{R}(M)\subseteq\overline{R}(F).$$
\end{dn}

 \begin{tm}\label{tm-FSAp-FSdom}  For an FS-approximation space $(U, R, \mathcal{F})$,  $(\mathfrak{C}(U, R, \mathcal{F}), \subseteq)$ is an FS-domain.
\end{tm}
 \noindent{\bf Proof.}   It follows from Lemma \ref{pn3-cf-cl}(2) and \ref{tm-cfga-way}(2) that the family $\{\overline{R}(F)\mid F\in \mathcal{F}\}$ is
a basis of $(\mathfrak{C}(U, R, \mathcal{F}), \subseteq)$ and that $(\mathfrak{C}(U, R, \mathcal{F}), \subseteq)$ is a  continuous domain. Next we prove that there is an approximate identity for $(\mathfrak{C}(U, R, \mathcal{F}), \subseteq)$ consisting of finitely separating maps. Let $\{\Theta_{i}\}_{i\in I}$  be  a directed family  of  CF-approximable relations on $(U, R, \mathcal{F})$ satisfying the condition (FS 1) and (FS 2) in Definition \ref{dn-FS-ga}. For all $i\in I$, the map $f_{\Theta_{i}}$ is a Scott continuous maps on $(\mathfrak{C}(U, R, \mathcal{F}), \subseteq)$ by Lemma \ref{tm-cfap-sco}. Let $\Theta_{j}, \Theta_{k}\in \{\Theta_{i}\}_{i\in I}$ and $\Theta_{j}\subseteq\Theta_{k}$. Notice that $f_{\Theta_{i}}(E)=\bigcup\{\overline{R}(G)\mid F, G\in \mathcal{F}, F\subseteq E\mbox{\ and\ } F\mathrel{\Theta}_i  G\}$ $(i\in I)$ for all $E \in\mathfrak{C}(U, R, \mathcal{F})$. Then we have $f_{\Theta_{j}}(E)\subseteq f_{\Theta_{k}}(E)$ for all $E \in\mathfrak{C}(U, R, \mathcal{F})$, showing that the  family $\{f_{\Theta_{i}}\}_{i\in I}$ is directed by the directedness of $\{\Theta_{i}\}_{i\in I}$. For all  $E\in\mathfrak{C}(U, R, \mathcal{F})$, we have
 \begin{align*}
(\bigvee_{i\in I} f_{\Theta_{i}})(E)=&\bigvee_{i\in I}(f_{\Theta_{i}}(E)) \mbox{ (by Lemma II-2.5 in }\cite{Gierz})\\
=&\bigcup_{i\in I}(f_{\Theta_{i}}(E))~~(\mbox{by~Lemma~}\ref{pn2-cf-clo}(3))\\
=&\bigcup_{i\in I}(\bigcup\{\overline{R}(G)\mid F, G\in \mathcal{F}, F\subseteq E \mbox{ and } F\mathrel{\Theta_{i}} G\})\\
=&\bigcup\{\overline{R}(G)\mid F, G\in \mathcal{F}, F\subseteq E \mbox{ and } (F, G)\in \bigcup_{i\in I}\Theta_{i} \}\\
=&\bigcup\{\overline{R}(G)\mid F, G\in \mathcal{F}, F\subseteq E \mbox{ and } (F, G)\in \operatorname{Id}_{(U, R, \mathcal{F})} \}\mbox{ (by (FS 1))}\\
=&\bigcup\{\overline{R}(G)\mid F, G\in \mathcal{F}, F\subseteq E \mbox{ and } G\subseteq\overline{R}(F) \}\\
=&\bigcup\{\overline{R}(F)\mid F\in \mathcal{F} \mbox{ and }F\subseteq E \}
(\mbox{by~Lemma~}\ref{pn2-cf-clo}(1)\mbox{ and }\ref{pn3-cf-cl}(2))\\
=& E \mbox{ (by Lemma }\ref{pn3-cf-cl}(2)\mbox{)}.
\end{align*}
This shows that $\bigvee_{i\in I} f_{\Theta_{i}}=id_{\mathfrak{C}(U, R, \mathcal{F})}$.

Next we verify that for all $i\in I$, $f_{\Theta_{i}}$ is  finitely separating. By the condition (FS 2) in Definition \ref{dn-FS-ga}, we have  a finite subfamily $\mathcal{M}_{i}$ of $\mathcal{F}$ such that for all $F\in \mathcal{F}$, there exists  $N\in \mathcal{M}_{i}$ satisfying
\vskip5pt
\centerline{$\forall G\in \mathcal{F}, F\mathrel{\Theta_{i}}G \Rightarrow G\subseteq \overline{R}(N)\subseteq\overline{R}(F).$}
 \vskip5pt
\noindent For $M\in \mathcal{M}_{i}$ and
 $E\in \mathfrak{C}(U, R, \mathcal{F})$, set
\vskip5pt
\centerline{$\mathcal{D}_{M}=\{\overline{R}(F)\mid  ( F\in \mathcal{F},  F\subseteq E) \mbox{ and } (\forall G\in \mathcal{F}, F\mathrel{\Theta_{i}} G\Rightarrow G\subseteq \overline{R}(M)\subseteq\overline{R}(F)\}.$}
\vskip5pt
\noindent Then we have $\bigcup_{M\in \mathcal{M}_{i}}\mathcal{D}_{M}=\{\overline{R}(F)\mid   F\in \mathcal{F} \mbox{ and }  F\subseteq E\}$. By Lemma \ref{pn3-cf-cl}(2),  family $ \{\overline{R}(F)\mid   F\in \mathcal{F} \mbox{ and }  F\subseteq E\}$ is directed.  Since $\mathcal{M}_{i}$ is finite, by Lemma \ref{lm-D=cup-direct}, there exists $M_{0}\in \mathcal{M}_{i}$ such that $\mathcal{D}_{M_{0}}$ is a cofinal subfamily of $\{\overline{R}(F)\mid   F\in \mathcal{F} \mbox{ and }  F\subseteq E\}$ and $\bigcup\mathcal{D}_{M_{0}}=\bigcup\{\overline{R}(F)\mid   F\in \mathcal{F} \mbox{ and }  F\subseteq E\}$. It follows from  Lemma \ref{pn3-cf-cl}(2) that $\bigcup\mathcal{D}_{M_{0}}= E$.
 Set $\overline{\mathcal{M}_{i}}=\{\overline{R}(M)\mid M\in\mathcal{M}_{i}\}$. For $F, G\in \mathcal{F}$ with $F\subseteq E \mbox{ and } F\mathrel{\Theta_{i}} G$, by finiteness of $F$, directedness of $\mathcal{D}_{M_{0}}$ and that $\bigcup\mathcal{D}_{M_{0}}=E$,  we can find some $F_{0}\in \mathcal{F}$ satisfying  $F_0\subseteq E$ and,
\vskip5pt
\centerline{$\forall G\in \mathcal{F}, F_0\mathrel{\Theta_{i}}G \Rightarrow G\subseteq \overline{R}(M_{0})\subseteq\overline{R}(F_0).$}
\vskip5pt
\noindent (hence $\overline{R}(F_{0})\in\mathcal{D}_{M_{0}}$) such that $F\subseteq  \overline{R}(F_{0})$.
 It follows from $F\subseteq  \overline{R}(F_{0})$ and  $F\mathrel{\Theta_{i}} G$ that $F_{0}\mathrel{\Theta_{i}} G$ by  Definition \ref{dn-CF-app}(2). Thus $G\subseteq \overline{R}(M_{0})\subseteq\overline{R}(F_0)$. It follows from   $G\subseteq\overline{R}(M_{0})$ that $\overline{R}(G)\subseteq\overline{R}(M_{0})$  by Lemma \ref{lm-ref-lm2-tr-up}(3).  Since
\vskip5pt
\centerline{$f_{\Theta_{i}}(E)=\bigcup\{\overline{R}(G)\mid F, G\in \mathcal{F}, F\subseteq E \mbox{ and } F\mathrel{\Theta_{i}} G\},$}
\vskip5pt
\noindent we have $f_{\Theta_{i}}(E)\subseteq\overline{R}(M_{0})\subseteq E$ and $\overline{R}(M_{0})\in \overline{\mathcal{M}_{i}}$. By finiteness of $\overline{\mathcal{M}_{i}}$, we see that $f_{\Theta_{i}}$  is  finitely separating.

 To sum up, it is proved that $(\mathfrak{C}(U, R, \mathcal{F}), \subseteq)$ is an FS-domain.\hfill$\Box$

\begin{dn}\label{dn-SFS-ga} A CF-approximation space $(U, R, \mathcal{F})$ is called {\em a strong FS-approximation space} if there exists a directed family $\{\Theta_{i}\}_{i\in I}$ of  CF-approximable relations from $(U, R, \mathcal{F})$ to itself satisfying the following conditions:

 (FS 1) $\bigcup_{i\in I}\Theta_{i}=\operatorname{Id}_{(U, R, \mathcal{F})}$;

(FS $2^{\prime}$)  For all $\Theta_{i}\ (i\in I)$, there is  $\mathcal{M}_{i}\subseteq_{fin}\mathcal{F}$ such that for all $F\in \mathcal{F}$, there exists  $M\in \mathcal{M}_{i}$ satisfying
$$\forall G\in \mathcal{F}, F\mathrel{\Theta_{i}}G \Rightarrow G\subseteq \overline{R}(M)\mbox{ and } M\subseteq\overline{R}(F).$$
\end{dn}

\begin{rk} \label{rk-SF-FS}
 By Lemma \ref{lm-ref-lm2-tr-up}(3), we have that (FS $2^{\prime}$)$\Rightarrow$(FS $2$). Thus a strong FS-approximation space $(U, R, \mathcal{F})$ must be an FS-approximation space, and $(\mathfrak{C}(U, R, \mathcal{F}), \subseteq)$ is an FS-domain. If $R$ is a preorder, then (FS $2^{\prime}$)$\Leftrightarrow$(FS $2$).
\end{rk}
Since for a preorder $R$ on a set $U$, the lower approximation operator $\underline{R}$  is really an interior operator for a topology on $U$, we have naturally the following definitions.
\begin{dn}\label{dn-cf-tga}   Let $(U, R)$ be a GA-space, $\emptyset\ne\mathcal{F}\subseteq \mathcal{P}_{fin}(U)\cup\{\emptyset\}$.  If $R$ is a preorder, then a CF-approximation space (resp., FS-approximation space) $(U, R, \mathcal{F})$ is called {\em a topological CF-approximation space} (resp., {\em topological FS-approximation space}).
\end{dn}
 \begin{rk}\label{GA-Pre-topCF} It follows from Lemma \ref{lm-ref-lm2-tr-up} that for a GA-space $(U, R)$ with $R$ being a preorder and $\emptyset\ne\mathcal{F}\subseteq \mathcal{P}_{fin}(U)\cup\{\emptyset\}$, then $(U, R, \mathcal{F})$ is a CF-approximation space, hence a topological CF-approximation space.
 \end{rk}

 \begin{ex}\label{ex-alg-top-cf-apps}
 Let $(L, \leqslant)$ be an algebraic  domain, $R_{K(L)}=\leqslant_{K(L)}$ and $\mathcal{F}_{K(L)}=\{F\subseteq_{fin} K(L)\mid F\  \mbox{has a top element } c_F\}$, where $\leqslant_{K(L)}$ is the restriction of $\leqslant$ to $K(L)$. Then

$(1)$ $(K(L), R_{K(L)},  \mathcal{F}_{K(L)})$ is a topological CF-approximation space;

$(2)$ $\mathfrak{C}(K(L), R_{K(L)},  \mathcal{F}_{K(L)})=\{{\downarrow} x\cap K(L)\mid x\in L\}$.
 \end{ex}
 \noindent {\bf Proof.} (1) Trivial.
 (2) Similar to the proof of Example \ref{tm-CDO-CFGA}(2).\hfill$\Box$\\

In this paper, we call $(K(L),  R_{K(L)}, \mathcal{F}_{K(L)})$ the induced topological CF-approximation space by the algebraic domain $L$.\\

The following proposition  gives a characterization of CF-approximable relation between topological CF-approximation spaces.
\begin{pn}\label{pn-TopCF-app}  Let $(U_{1}, R_{1}, \mathcal{F}_{1})$, $(U_{2}, R_{2}, \mathcal{F}_{2})$ be topological CF-approximation spaces. Then  a binary relation $\mathrel{\Theta}\subseteq\mathcal{F}_{1}\times\mathcal{F}_{2}$ is a CF-approximable relation iff it satisfies

$(1)$  for all $F\in \mathcal{F}_1$, there is $G\in \mathcal{F}_{2}$ such that $F\mathrel{\Theta} G$;

$(2)$  for all $F, F^{\prime}\in \mathcal{F}_{1}$, $G, G^{\prime}\in  \mathcal{F}_{2}$, if $F\subseteq \overline{R_{1}}(F^{\prime})$, $G^{\prime}\subseteq \overline{R_{2}}(G)$, and $F\mathrel{\Theta} G$, then $F^{\prime}\mathrel{\Theta} G^{\prime}$; and

$(3)$  for all $F\in \mathcal{F}_{1}$, $G_{1}, G_{2}\in \mathcal{F}_{2}$, if $F\mathrel{\Theta} G_{1}$ and $F\mathrel{\Theta} G_{2}$, then there is  $G_{3}\in \mathcal{F}_{2}$ such that  $G_{1}\cup G_{2}\subseteq\overline{R_{2}}(G_{3})$ and $F\mathrel{\Theta} G_{3}$.
\end{pn}
\noindent{\bf Proof.} Follows directly from Lemma \ref{lm-ref-lm2-tr-up} and that $R$ is a preorder.\hfill$\Box$
\begin{tm}\label{tm-tpf-agl}
If $(U, R, \mathcal{F})$ is a topological  FS-approximation space, then $(\mathfrak{C}(U, R, \mathcal{F}), \subseteq)$ is a BF-domain.
\end{tm}
\noindent{\bf Proof.} It follows from Lemma \ref{lm-ref-lm2-tr-up}(1), \ref{pn3-cf-cl}(2)  and \ref{tm-cfga-way}(1) that the family $\{\overline{R}(F)\mid F\in \mathcal{F}\}$ is
a basis of $(\mathfrak{C}(U, R, \mathcal{F}), \subseteq)$ consisting of compact elements. Thus $(\mathfrak{C}(U, R, \mathcal{F}), \subseteq)$ is an algebraic domain. By Theorem \ref{tm-FSAp-FSdom}, we know that $(\mathfrak{C}(U, R, \mathcal{F}), \subseteq)$ is an FS-domain,  hence a  BF-domain.
 \hfill$\Box$

\section{Representations of FS-domains and BF-domains}
In this section, we first consider  representations of FS-domains via  FS-approximation spaces. Then we introduce  topological  BF-approximation spaces and consider  representations of BF-domains via topological  BF-approximation spaces.

%

\begin{pn}\label{pn-dom-map}
Let $(L_{1}, R_{L_{1}},  \mathcal{F}_{L_{1}})$ and $(L_{2}, R_{L_{2}}, \mathcal{F}_{L_{2}})$ be the  induced CF-approximation spaces by continuous domains $L_{1}$ and $L_{2}$, respectively. For a  Scott continuous map $g: L_{1}\longrightarrow L_{2}$, define  $\Omega_{g}\subseteq\mathcal{F}_{L_{1}}\times\mathcal{F}_{L_{2}}$ such that
 \vskip3pt

 \centerline{$\forall F\in \mathcal{F}_{L_{1}},   \forall G\in \mathcal{F}_{L_{2}}, (F, G)\in \Omega_{g}\Leftrightarrow c_{G}\ll_{L_{2}} g(c_{F}), $}
 \vskip5pt

\noindent where $c_F$ is the top element of $F$. Then we have

 $(1)$ $\Omega_{g}$ is a CF-approximable relation from $(L_{1}, R_{L_{1}}, \mathcal{F}_{L_{1}})$ to $(L_{2}, R_{L_{2}}, \mathcal{F}_{L_{2}})$; and

   $(2)$ if $h\in [L_{1}\longrightarrow L_{2}]$ and $g\leqslant h$, then   $\Omega_{g}\subseteq \Omega_{h}$, that is if $(F, G)\in\Omega_{g}$, then  $(F, G)\in\Omega_{h}$.
\end{pn}
\noindent{\bf Proof.}
  (1) We need to prove that $\Omega_{g}$  satisfies the conditions (1)-(5) in Definition \ref{dn-CF-app}.

 Let $F\in \mathcal{F}_{L_{1}}$. Since $L_{2}$ is a continuous domain, we know that $\dda g(c_{F})\neq\emptyset$. Therefore there exists $x\in L_{2}$ such that $x\ll g(c_{F})$. Thus $\{x\}\in \mathcal{F}_{L_{2}}$ and $(F, \{x\})\in \Omega_{g}$.  This shows that $\Omega_{g}$  satisfies the condition (1) in Definition \ref{dn-CF-app}.

  To check that $\Omega_{g}$ satisfies the condition (2) in Definition \ref{dn-CF-app}, let $F, F^{\prime}\in \mathcal{F}_{L_{1}}$ and  $G\in \mathcal{F}_{L_{2}}$. Then
\begin{align*}
 &F\subseteq \overline{R_{L_{1}}}(F^{\prime})=\dda c_{F^{\prime}}, (F, G)\in \Omega_{g}\\
\Rightarrow& c_{F}\ll c_{F^{\prime}},  c_{G}\ll g(c_{F}) \\
\Rightarrow& c_{G}\ll g(c_{F})\leqslant g(c_{F^{\prime}})\\
\Rightarrow&(F^{\prime}, G)\in \Omega_{g}.
\end{align*}

To check that $\Omega_{g}$ satisfies the condition (3) in Definition \ref{dn-CF-app},  let $F\in \mathcal{F}_{L_{1}}$ and $G, G^{\prime}\in \mathcal{F}_{L_{2}}$. Then
\begin{align*}
&(F, G)\in \Omega_{g}, G^{\prime}\subseteq\overline{R_{L_{2}}}(G)=\dda c_{G} \\
\Rightarrow& c_{G^{\prime}}\ll c_{G},  c_{G}\ll g(c_{F}) \\
\Rightarrow& c_{G^{\prime}}\ll g(c_{F})\\
\Rightarrow&(F, G^{\prime} )\in \Omega_{g}.
\end{align*}

To check that $\Omega_{g}$ satisfies the  condition (4) in Definition \ref{dn-CF-app},  let $F\in \mathcal{F}_{L_{1}}$, $G\in \mathcal{F}_{L_{2}}$ and $(F, G)\in \Omega_{g}$.  Then  $c_{G}\ll g(c_{F})$. By Lemma \ref{pn-int}, there exist $a_{1}, a_{2}\in L_{2}$ such that $c_{G}\ll {a_{1}}\ll {a_{2}}\ll g(c_{F})$. Since $L_{1}$ is a continuous domain and $g$ is Scott continuous, we have that $g( c_{F})
=\bigvee\{g(x)\mid x\ll c_{F}\}$ and $\{g(x)\mid x\ll c_{F}\}$ is directed. It follows from  ${a_{2}}\ll g(c_{F})=\bigvee\{g(x)\mid x\ll c_{F}\}$ that there exists $v\ll c_{F}$ such that $a_{2}\leqslant g(v)$. Thus $c_{G}\ll a_{1}\ll g(v)$. Set $F^{\prime}=\{v\}\in  \mathcal{F}_{L_{1}}$ and $G^{\prime}=\{a_{1}\}\in  \mathcal{F}_{L_{2}}$. Then $G\subseteq \dda a_{1} =\overline{R_{L_{2}}}(G^{\prime})$, $F^{\prime}=\{v\}\subseteq\dda c_{F}=\overline{R_{L_{1}}}(F)$ and $(F^{\prime}, G^{\prime})\in \Omega_{g}$, showing that  $\Omega_{g}$ satisfies the condition (4).

To check that $\Omega_{g}$ satisfies  the condition (5) in Definition \ref{dn-CF-app}, let $F\in \mathcal{F}_{L_{1}}$ and  $G_{1}, G_{2}\in \mathcal{F}_{L_{2}}$ satisfying $(F, G_{1})\in \Omega_{g}$ and $(F, G_{2})\in \Omega_{g}$. Then $\{c_{G_{1}}, c_{G_{2}}\}\subseteq \dda g(c_{F})$. Since $L$ is a continuous domain, we have $d\in L_{2}$ such that $\{c_{G_{1}}, c_{G_{2}}\}\subseteq \dda d$ and $d\ll g(c_{F})$. Therefore $G_{1}\cup G_{2}\subseteq\overline{R_{L_{2}}}(\{d\})$ and $(F, \{d\})\in \Omega_{g}$, showing that  $\Omega_{g}$ satisfies the condition (5).

  (2) Straightforward.
\hfill$\Box$

%
%

\begin{tm}\label{tm-FSDO-FSGA} Let $(L, \leqslant)$ be an FS-domain. Then the induced CF-approximation space $(L, R_{L},  \mathcal{F}_{L})$ is an FS-approximation space.
\end{tm}
\noindent{\bf Proof.} Let $\{\delta_{i}\}_{i\in I}$ be an approximate identity for $L$ consisting of finitely separating maps. For all $i\in I$,  construct a binary relation $\Theta_{i}= \Omega_{\delta_{i}}$ such that
\vskip5pt
 \centerline{$ \forall F, G\in \mathcal{F}_{L},   (F, G)\in \Theta_{i}\Leftrightarrow c_{G}\ll \delta_{i}(c_{F}).$}
\vskip5pt
 \noindent It follows from Proposition \ref{pn-dom-map} that $\{\Theta_{i}\}_{i\in I}$ is a directed family of CF-approximable relations from $(L, R_{L},  \mathcal{F}_{L})$  to itself.  Next, we check that $\{\Theta_{i}\}_{i\in I}$ satisfies (FS 1) and (FS 2) in Definition \ref{dn-FS-ga}.

To check that $\{\Theta_{i}\}_{i\in I}$ satisfies (FS 1) in Definition \ref{dn-FS-ga}, let $F, G\in \mathcal{F}_{L}$ and
$(F, G)\in \operatorname{Id}_{(L, R_{L},  \mathcal{F}_{L}) }$. Then $G\subseteq\overline{R_{L}}(F)=\dda c_F$ and
$c_{G}\ll c_{F}$.
Since $\{\delta_{i}\}_{i\in I}$ is an approximate identity for $L$, we have that $\bigvee_{i\in I}( \delta_{i}(c_{F}))=c_{F}$.  Since $L$ is continuous and $c_{G}\ll c_{F}=\bigvee_{i\in I}( \delta_{i}(c_{F}))$, by interpolation property, there exists $i\in I$ such that $c_{G}\ll \delta_{i}(c_{F})$ and $(F, G)\in \Theta_{i}=\Omega_{\delta_{i}}$,  showing that $\operatorname{Id}_{(L, R_{L},  \mathcal{F}_{L})}\subseteq \bigcup_{i\in I} \Theta_{i}$.

Conversely, let $(F, G)\in \bigcup_{i\in I} \Theta_{i}$. Then there exists $j\in I$ such that $(F, G)\in \Theta_{j}=\Omega_{\delta_{j}}$. Therefore $c_{G}\ll \delta_{j}(c_{F})\leqslant c_{F}$. Thus $G\subseteq\dda c_{F}=\overline{R_{L}}(F)$. It follows from the definition of $\operatorname{Id}_{(L, R_{L},  \mathcal{F}_{L})}$ that $(F, G)\in  \operatorname{Id}_{(L, R_{L},  \mathcal{F}_{L})}$, showing that $\bigcup_{i\in I} \Theta_{i}\subseteq \operatorname{Id}_{(L, R_{L},  \mathcal{F}_{L})}$. Thus $\bigcup_{i\in I} \Theta_{i}= \operatorname{Id}_{(L, R_{L},  \mathcal{F}_{L})}$.  Thus (FS 1) is checked.

Next we  check that $\{\Theta_{i}\}_{i\in I}$ satisfies (FS 2) in Definition \ref{dn-FS-ga}. By Definition \ref{dn-doms}, for every $i\in I$, we can get  a finite set $M_{i}$ such that  for each $x\in L$, there exists $m\in M_{i}$ such that $\delta(x)\leqslant m\leqslant  x$. Construct a finite family $\mathcal{M}_{i}=\{\{m\}\mid m\in M_{i}\}$ for all $i\in I$. Since $\delta_{i}$ is finitely separating with $M_i$, we know that for every $F\in \mathcal{F}_{L}$, there exists $m\in M_{i}$ such that $\delta_{i}(c_{F})\leqslant m\leqslant c_{F}$. Now,  set $M=\{m\}\in \mathcal{M}_{i}$. If $G\in \mathcal{F}_{L}$ and $F\mathrel{\Theta_{i}} G$, then $c_{G}\ll \delta_{i}(c_{F})$.   Therefore $c_{G}\ll m\leqslant c_{F}$ and $G\subseteq \overline{R_{L}}(M)\subseteq \overline{R_{L}}(F)$, showing that $\Theta_{i}$ satisfies (FS 2) in Definition \ref{dn-FS-ga} for every $i\in I$.
\hfill$\Box$



\begin{tm}\label{tm-rep-FS} {\em (Representation Theorem I: for FS-domain)} A poset $(L, \leqslant)$ is an FS-domain iff there is an FS-approximation space $(U, R, \mathcal{F})$ such that $(L, \leqslant)\cong (\mathfrak{C}(U, R, \mathcal{F}), \subseteq))$.
\end{tm}
\noindent{\bf Proof.} $\Leftarrow$: Follows directly from Theorem \ref{tm-FSAp-FSdom}.

$\Rightarrow$: If $L$ is an FS-domain, then by Theorem \ref{tm-FSDO-FSGA},  the induced CF-approximation space $(L, R_{L},  \mathcal{F}_{L})$ is an FS-approximation space. Define  $f: L\to \mathfrak{C}(L, R_{L},  \mathcal{F}_{L})$ such that for all $x\in L$, $f(x)=\dda x$. Then it follows from  Example \ref{tm-CDO-CFGA}(2) and the continuity of $L$ that $f$ is an order
isomorphism.
\hfill$\Box$

\begin{tm}\label{tm-rep-FS2} {\em (Representation Theorem II: for FS-domain)} A poset $(L, \leqslant)$ is an FS-domain  iff there is a strong FS-approximation space $(U, R, \mathcal{F})$ such that $(L, \leqslant)\cong (\mathfrak{C}(U, R, \mathcal{F}), \subseteq))$.
\end{tm}
\noindent{\bf Proof.} $\Leftarrow$: Follows directly from Theorem \ref{tm-rep-FS} and Remark \ref{rk-SF-FS}.

$\Rightarrow$: By the proof of Theorem \ref{tm-rep-FS}
we need only to prove that $(L, R_{L},  \mathcal{F}_{L})$ defined in Example \ref{tm-CDO-CFGA} is a strong FS-approximation space, Let $\{\delta_{i}\}_{i\in I}$ be an approximate identity for $L$ consisting of finitely separating maps. For all $i\in I$,  construct a binary relation $\Theta_{i}= \Omega_{\delta_{i}^{2}}$
such that\\
 \centerline{$ \forall F, G\in \mathcal{F}_{L},   (F, G)\in \Theta_{i}\Leftrightarrow c_{G}\ll \delta_{i}^{2}(c_{F}).$}\\
 By Proposition \ref{pn-dom-map}, $\{\Theta_{i}\}_{i\in I}$ is a directed family of CF-approximable relations on $(L, R_{L},  \mathcal{F}_{L})$. 

Similar to the proof of Theorem \ref{tm-FSDO-FSGA}, it is easy to check that $\{\Theta_{i}\}_{i\in I}$ satisfies (FS 1) in Definition \ref{dn-SFS-ga} by Lemma \ref{lm-fis-map}(1).
Next we  check that $\{\Theta_{i}\}_{i\in I}$ satisfies (FS $2^{\prime}$) in Definition \ref{dn-SFS-ga}. By Definition \ref{dn-doms}, for every $i\in I$, we can get  a finite set $M_{i}$ such that  for each $x\in L$, there exists $m\in M_{i}$ such that $\delta(x)\leqslant m\leqslant  x$. Construct a finite family $\mathcal{M}_{i}=\{\{\delta_{i}(m)\}\mid m\in M_{i}\}$ for every $i\in I$. Since $\delta_{i}$ is finitely separating with $M_{i}$, we know that for every $F\in \mathcal{F}_{L}$, there exists $m\in M_{i}$ such that $\delta_{i}(c_{F})\leqslant m\leqslant c_{F}$. Now, set $M=\{\delta_{i}(m)\}\in \mathcal{M}_{i}$. If $G\in \mathcal{F}_{L}$ and $F\mathrel{\Theta}_{i} G$, then $c_{G}\ll \delta_{i}^{2}(c_{F})$. Noticing that $\delta_{i}(c_{F})\leqslant m\leqslant c_{F}$ and $\delta_{i}$ is order preserving, we have that $ \delta_{i}^{2}(c_{F})\leqslant  \delta_{i}(m)\leqslant\delta_{i}(c_{F})\ll c_{F}$ by Lemma \ref{lm-fis-map}. Thus $c_{G}\ll \delta_{i}(m)\ll c_{F}$, $G\subseteq \overline{R_{L}}(M)$ and $M\subseteq \overline{R_{L}}(F)$, showing that $\Theta_{i}$ satisfies (FS $2^{\prime}$) in Definition \ref{dn-SFS-ga} for every $i\in I$.
\hfill$\Box$

\begin{dn} (\cite{Gierz})\label{dn-ker-map}
Let $P$ be a poset,  $k: P\longrightarrow P$  a monotone map. If $k$ satisfies

 (1) $k(x)\leqslant x$ for all $x\in P$; and

 (2)  $k(k(x))=k(x)$ for all  $x\in P$,

\noindent then $k$ is called {\em a kernel operator}.
\end{dn}
\begin{lm}{\em (\cite[Proposition II-2.20]{Gierz})}\label{lm-BF-eq}
For a dcpo  $L$, the following properties are equivalent:

$(1)$ $L$ is a BF-domain;

$(2)$ $L$ is an algebraic domain and has an approximate identity consisting of  maps with finite range;

 $(3)$ $L$ has an approximate identity consisting of  kernel operators with finite range.
\end{lm}

Now we arrive at representations of BF-domain with topological FS-approximation spaces.

\begin{tm}  {\em (Representation Theorem III: for BF-domains)}\label{rep-BF-topF}
A poset $(L, \leqslant)$ is  a BF-domain iff  there is a topological FS-approximation space $(U, R, \mathcal{F})$ such that $(\mathfrak{C}(U, R, \mathcal{F}), \subseteq)\cong (L, \leqslant)$.
\end{tm}
\noindent{\bf Proof.} $\Leftarrow$:   Follows from Theorem \ref{tm-tpf-agl}.

$\Rightarrow$:  Let $(L, \leqslant)$  a BF-domain, $\mathcal{F}_{K(L)}=\{F\subseteq_{fin} K(L)\mid F\ \mbox{has a top element}\}$ and $R_{K(L)}=\leqslant_{K(L)}$. By Remark \ref{GA-Pre-topCF}, $(K(L),  R_{K(L)}, \mathcal{F}_{K(L)})$ is a topological CF-approximation space. Moreover, for any $F\in \mathcal{F}_{K(L)}$, let $c_{F}$ be the top element of $F$. To show that $(K(L),  R_{K(L)}, \mathcal{F}_{K(L)})$ is a topological FS-approximation space, 
we need only to prove  $(K(L),  R_{K(L)}, \mathcal{F}_{K(L)})$ satisfies (FS 1) and (FS 2). By Lemma \ref{lm-BF-eq}, we have  an approximate identity $\{\delta_{i}\}_{i\in I}$  for $L$ consisting of kernel operators with finite range. For all $i\in I$, define $\Theta_{i}\subseteq\mathcal{F}_{K(L)}\times\mathcal{F}_{K(L)}$ such that
$$\forall F, G\in\mathcal{F}_{K(L)}, (F, G)\in \Theta_{i}\Longleftrightarrow c_{G}\leqslant \delta_{i}(c_{F}).$$
It is routine to check that $\{\Theta_{i}\}_{i\in I}$ is a directed family of CF-approximable relations.

 To check that $\{\Theta_{i}\}_{i\in I}$  satisfies (FS 1), let $F, G\in  \mathcal{F}_{K(L)}$, then
\begin{align*}
(F, G)\in \bigcup_{i\in I}\Theta_{i}\Longleftrightarrow&\exists j\in I \mbox{ such that }(F, G)\in \Theta_{j}\\
\Longleftrightarrow& \exists j\in I\mbox{ such that }c_{G}\leqslant \delta_{j}(c_{F})\\
\Longleftrightarrow& c_{G}\leqslant c_{F}(``\Leftarrow" \mbox{ follows by }\bigvee_{i\in I}\delta_{i}(c_F)=c_{F} \mbox{ and } c_{G}\in K(L))\\
\Longleftrightarrow& G\subseteq\overline{R_{K(L)}}(F)={\downarrow} c_{F}\cap K(L)\\
\Longleftrightarrow& (F, G)\in \operatorname{Id}_{(K(L),  R_{K(L)}, \mathcal{F}_{K(L)})}.
\end{align*}

Next we  check that $\{\Theta_{i}\}_{i\in I}$  satisfies (FS 2). For all $\Theta_{i}$, let $M_{i}$ be the range of $\delta_{i}$. For any $m\in M_i$, by Lemma \ref{lm-fis-map}(2) and the condition (2) in Definition \ref{dn-ker-map}, we have that $m=\delta_{i}(m)\ll m$ and $m$ is compact. Set $\mathcal{M}_{i}=\{\{m\}\mid m\in M_{i}\}\subseteq \mathcal{F}_{K(L)}$. For all $F\in  \mathcal{F}_{K(L)}$,  select $\{\delta_{i}(c_F)\}\in \mathcal{M}_{i}$. If $G\in  \mathcal{F}_{K(L)}$ and $F\mathrel{\Theta_{i}} G$, then $c_{G}\leqslant\delta_{i}(c_F)\leqslant c_F$, showing that $G\subseteq \overline{R_{\small{K(L)}}}(\{\delta_{i}(c_F)\})\subseteq \overline{R_{K(L)}}(F)$. Thus $\{\Theta_{i}\}_{i\in I}$  satisfies (FS 2).

  By Example \ref{ex-alg-top-cf-apps}(2), we see that $(\mathfrak{C}(K(L),  R_{K(L)}, \mathcal{F}_{K(L)}), \subseteq)\cong (L, \leqslant)$.
\hfill$\Box$\\

 Representations of FS-domains and BF-domains above are all given by using CF-approximable relations. In the following, we will give direct representations of  BF-domains without using CF-approximable relations. To this end, we give
 another one of main concepts in this paper as follows.
\begin{dn}\label{dn-top-BF} {\em A topological BF-approximation space} is  a topological CF-approximation space $(U, R, \mathcal{F})$ satisfying that  for all $K\subseteq_{fin}U$, there  is a finite family $\mathcal{M}_{K}\subseteq_{fin}\mathcal{F}$ such that

(TB 1) $\mathcal{P}(K)\cap \mathcal{F}\subseteq \mathcal{M}_{K}$;

 (TB 2) $(\forall F\in \mathcal{F}) (\mathcal{G}\subseteq\mathcal{M}_{K}, \bigcup\mathcal{G}\subseteq \overline{R}(F))\Longrightarrow (\exists M\in \mathcal{M}_{K}) (\bigcup\mathcal{G}\subseteq \overline{R}(M)\subseteq \overline{R}(F))$.
\end{dn}
\begin{rk}\label{rk-mk-bfs} For a topological BF-approximation space $(U, R, \mathcal{F})$, we have

 $(1)$ the binary relation $R$ is a preorder;

 $(2)$ for the same  $K\subseteq_{fin}U$, we have  may more than one   finite families of $\mathcal{F}$ satisfying  conditions (TB 1) and (TB 2), however,  we can use Axiom of Choice for all $K\subseteq_{fin}U$ to select a fixed one $\mathcal{M}_{K}\subseteq_{fin}\mathcal{F}$ satisfying  conditions (TB 1) and (TB 2).

 $(3)$ for all $K\subseteq_{fin}U$, set $\mathcal{G}=\emptyset$, then by  (TB 2) in Definition \ref{dn-top-BF}, we have that $\mathcal{M}_{K}\neq\emptyset$.
\end{rk}
\begin{tm}\label{tm-ToBF-BF}  If $(U, R, \mathcal{F})$ is a topological BF-approximation space,  then $(\mathfrak{C}(U, R, \mathcal{F}), \subseteq)$ is a BF-domain.
\end{tm}
\noindent{\bf Proof.} Clearly, $(\mathfrak{C}(U, R, \mathcal{F}), \subseteq)$ is an algebraic domain. To show that $(\mathfrak{C}(U, R, \mathcal{F}), \subseteq)$ is a BF-domain, we divide the  proof into several steps.

Step 1. Set $ \mathcal{D}=\{K\subseteq_{fin} U\mid \exists F\in \mathcal{F}\mbox{~s.t.~} F\subseteq K\}$. Then in  set-inclusion order, $ \mathcal{D}$ is clearly a directed family.
For all $K\in \mathcal{D}$, define
$\delta_{K}: \mathfrak{C}(U, R, \mathcal{F}) \longrightarrow \mathfrak{C}(U, R, \mathcal{F})$ such that for all $E\in \mathfrak{C}(U, R, \mathcal{F})$, $\delta_{K}(E)=\bigcup\{\overline{R}(M)\mid M\in \mathcal{M}_{K}\mbox{ and } M\subseteq E\}$, where $\mathcal{M}_{K}$ is stated in Remark \ref{rk-mk-bfs}(2).

Step 2. Assert that for all $K\in \mathcal{D}$, $\delta_{K}$ is well defined and has  finite range.

Let $E\in  \mathfrak{C}(U, R, \mathcal{F})$. By the finiteness of $\mathcal{M}_{K}$ and $\mathcal{M}_{K}\subseteq_{fin}\mathcal{F}$, we have $\bigcup\{M\mid M\in \mathcal{M}_{K}\mbox{ and } M\subseteq E\}\subseteq_{fin} E$. It follows from $E\in \mathfrak{C}(U, R, \mathcal{F})$ that there exists $F\in \mathcal{F}$ such that $\bigcup\{M\mid M\in \mathcal{M}_{K}\mbox{ and } M\subseteq E\}\subseteq \overline{R}(F)\subseteq E$.   It follows from $\{M\mid M\in \mathcal{M}_{K}\mbox{ and } M\subseteq E\}\subseteq \mathcal{M}_{K}$ and (TB 2) that there is $M^{*}\in \mathcal{M}_{K}$ such that $\bigcup\{M\mid M\in \mathcal{M}_{K}\mbox{ and } M\subseteq E\}\subseteq \overline{R}(M^{*})\subseteq\overline{R}(F)\subseteq E$. By Lemma \ref{lm-ref-lm2-tr-up} (1) and (3), we see that  $M^{*}\subseteq E$ and $\overline{R}(M^{*})$ is the greatest element in $\{\overline{R}(M)\mid M\in \mathcal{M}_{K}\mbox{ and } M\subseteq E\}$ equipped with set-inclusion order. Hence $\delta_{K}(E)=\overline{R}(M^{*})\in \mathfrak{C}(U, R, \mathcal{F})$, showing that $\delta_{K}$ is well defined. It follows from the finiteness of $\mathcal{M}_{K}$ and $M^{*}\in \mathcal{M}_{K}$ that $\delta_{K}$ has  finite range.

Step 3.  Assert that for all $K\in \mathcal{D}$, $\delta_{K}$ is Scott continuous.

Obviously, $\delta_{K}$ is order preserving. Let $\{E_{i}\}_{i\in I}\subseteq(\mathfrak{C}(U, R, \mathcal{F}), \subseteq)$ be a  directed family and  $E=\bigvee_{i\in I} E_{i}=\bigcup_{i\in I}E_{i}\in\mathfrak{C}(U, R, \mathcal{F})$. By the proof of Step 2, there is $M_{E}^{*}\in \mathcal{M}_{K}$ and $M_{E}^{*}\subseteq E$ such that $\delta_{K}(E)=\overline{R}(M_{E}^{*})$. Since $M_{E}^{*}\subseteq_{fin}\bigcup_{i\in I}E_{i}=E$, there exists $j\in I$ such that  $M_{E}^{*}\subseteq E_{j}$.  So,   $\delta_{K}(E_{j})=\bigcup\{\overline{R}(M)\mid M\in \mathcal{M}_{K}\mbox{ and } M\subseteq E_{j}\}\supseteq \overline{R}(M_{E}^{*})=\delta_{K}(E)$. Therefore $\delta_{K}(E)\subseteq\delta_{K}(E_{j})\subseteq\bigcup_{i\in I}\delta_{K}(E_{i})\subseteq\delta_{K}(E)$. Thus $\delta_{K}(\bigcup_{i\in I}E_{i})=\bigcup_{i\in I}\delta_{K}(E_{i})$, showing that $\delta_{K}$ is Scott continuous.

Step 4. Assert that $\{\delta_{K}\}_{K\in \mathcal{D}}$ is an approximate identity on $(\mathfrak{C}(U, R, \mathcal{F}), \subseteq)$.

To show $\{\delta_{K}\}_{K\in \mathcal{D}}$ is directed, let $K_{1}, K_{2}\in \mathcal{D}$ and $K=K_{1}\cup K_{2}\cup\bigcup(\mathcal{M}_{K_{1}}\cup\mathcal{M}_{K_{2}})$. Then $K\in \mathcal{D}$ and $\mathcal{M}_{K_{1}}\cup\mathcal{M}_{K_{2}}\subseteq \mathcal{P}(K)\cap \mathcal{F}\subseteq\mathcal{M}_{K}$. For all $E\in \mathfrak{C}(U, R, \mathcal{F})$, $i=1, 2$,  we have that
\vskip 3pt
\centerline{$\delta_{K_{i}}(E)=\bigcup\{\overline{R}(M)\mid M\in \mathcal{M}_{K_{i}}\mbox{ and } M\subseteq E\}\subseteq\bigcup\{\overline{R}(M)\mid M\in \mathcal{M}_{K}\mbox{ and } M\subseteq E\}=\delta_{K}(E),$}
\vskip 3pt
\noindent showing that $\{\delta_{K}\}_{K\in \mathcal{D}}$ is directed.

Let $F\in \mathcal{F}$ and $F\subseteq E\in \mathfrak{C}(U, R, \mathcal{F})$. Then by (TB 1), there exists $\mathcal{M}_{F}\subseteq_{fin}\mathcal{F}$ satisfying $\mathcal{P}(F)\cap \mathcal{F}\subseteq \mathcal{M}_{F}$. It follows from  $F\in \mathcal{P}(F)\cap \mathcal{F}\subseteq \mathcal{M}_{F}$  that $\overline{R}(F)\subseteq\bigcup\{\overline{R}(M)\mid M\in \mathcal{M}_{F}\mbox{ and } M\subseteq E\}=\delta_{F}(E)$. Noticing that $F\in \mathcal{F}\subseteq\mathcal{D}$, we have $\overline{R}(F)\subseteq\bigvee_{K\in \mathcal{D}}\delta_{K}(E)$. By Lemma \ref{pn3-cf-cl} and the arbitrariness $F$, we have  $E=\{\overline{R}(F)\mid F\in \mathcal{F}, F\subseteq E\}\subseteq\bigvee_{K\in \mathcal{D}}\delta_{K}(E)$. Obviously, $\bigvee_{K\in \mathcal{D}}\delta_{K}(E)\subseteq E$ and $\bigvee_{K\in \mathcal{D}}\delta_{K}(E)= E$. This shows that $\bigvee_{K\in \mathcal{D}}\delta_{K}=id_{\mathfrak{C}(U, R, \mathcal{F})}$.

To sum up Step 1--Step 4, by Lemma \ref{lm-BF-eq}, we have that $(\mathfrak{C}(U, R, \mathcal{F}), \subseteq)$ is a BF-domain.\hfill$\Box$\\

Since a BF-domain is an FS-domain, it is natural to ask that whether  a topological  BF-approxiamtion space is a topological  FS-approxiamtion space or not? The following proposition gives an affirmative answer.
\begin{pn}\label{pn-topBF>topFS}  A topological BF-approxiamtion space  is  a topological FS-approxiamtion space.
\end{pn}
\noindent{\bf Proof.} Let $(U, R, \mathcal{F})$ be a topological BF-approxiamtion space. Then  $\mathcal{D}=\{K\subseteq_{fin} U\mid \exists F\in \mathcal{F}\mbox{ such that } F\subseteq K\}$ defined in the proof of Theorem \ref{tm-ToBF-BF} is a directed family.  For all $K\in \mathcal{D}$,  define a binary  relation $\Theta_{K}$ on $ \mathcal{F}$ such that

\centerline{ $\forall F, G\in \mathcal{F}, (F, G)\in \Theta_{K}\Longleftrightarrow \exists M\in \mathcal{M}_{K}\mbox{ such that } G\subseteq \overline{R}(M)\subseteq \overline{R}(F),$}
 \noindent where $\mathcal{M}_{K}$ is stated in Remark \ref{rk-mk-bfs}(2).

To show the proposition, it suffices by Definition \ref{dn-FS-ga} to  show that the family $\{\Theta_{K}\}_{K\in \mathcal{D}}$ is a directed family of  CF-approximable relations on $(U, R, \mathcal{F})$ satisfying  (FS 1) and (FS 2) in Definition \ref{dn-FS-ga}. To this end, we divide the proof into several steps.

Step 1.  Show that for all $K\in \mathcal{D}$, $\Theta_{K}$ is a CF-approximable relation.

To check $\Theta_{K}$ satisfies the condition (1) in Proposition \ref{pn-TopCF-app}, let $F\in \mathcal{F}$. By the proof of Theorem \ref{tm-ToBF-BF},  there is an $M^{*}\in \mathcal{M}_{K}$ such that $\overline{R}(M^{*})$ is the greatest element in $\{\overline{R}(M)\mid M\in \mathcal{M}_{K}\mbox{ and } M\subseteq \overline{R}(F)\}$. Thus $M^{*}\subseteq\overline{R}(M^{*})\subseteq \overline{R}(F)$ and $(F, M^{*})\in \Theta_{K}$, showing  $\Theta_{K}$ satisfies the condition (1).

To check $\Theta_{K}$ satisfies the condition (2) in Proposition \ref{pn-TopCF-app}, let  $F, F^{\prime}\in \mathcal{F}_{1}$, $G, G^{\prime}\in  \mathcal{F}_{2}$. Then
 \begin{align*}
 & F\subseteq \overline{R}(F^{\prime}), G^{\prime}\subseteq \overline{R}(G) \mbox{ and } F\mathrel{\Theta_{K}} G\\
\Longrightarrow& \exists M\in \mathcal{M}_{K}, F\subseteq \overline{R}(F^{\prime}), G^{\prime}\subseteq \overline{R}(G) \mbox{ and }G\subseteq \overline{R}(M)\subseteq \overline{R}(F)\\
\Longrightarrow& \exists M\in \mathcal{M}_{K}, G^{\prime}\subseteq \overline{R}(M)\subseteq \overline{R}(F^{\prime})\\
\Longrightarrow& F^{\prime}\mathrel{\Theta_{K}} G^{\prime}.
\end{align*}
So, $\Theta_{K}$ satisfies the condition (2).

To check $\Theta_{K}$ satisfies the condition (3) in Proposition \ref{pn-TopCF-app},  let $F\in \mathcal{F}_{1}$, $G_{1}, G_{2}\in \mathcal{F}_{2}$.
If $(F, G_{1})\in\mathrel{\Theta_{K}}$ and $(F, G_{2})\in \mathrel{\Theta_{K}}$,  then there exists $M_1, M_2\in \mathcal{M}_{K}$ such that $G_i\subseteq \overline{R}(M_i)\subseteq \overline{R}(F)$ $(i=1, 2)$. By the maximality of $\overline{R}(M^*)$, we have that $G_1\cup G_2\subseteq \overline{R}(M_1)\cup \overline{R}(M_2)\subseteq \overline{R}(M^*)\subseteq\overline{R}(F)$ and $(F, M^*)\in\mathrel{\Theta_{K}}$, showing that  $\Theta_{K}$ satisfies the condition (3).

Thus, $\Theta_{K}$ is a CF-approximable relation on $(U, R, \mathcal{F})$.

Step 2. It follows from Step 4 in the proof of Theorem \ref{tm-ToBF-BF} that $(\{\mathcal{M}_{K}\}_{K\in \mathcal{D}}, \subseteq)$ is directed. Thus $\{\Theta_{K}\}_{K\in \mathcal{D}}$ is directed.

Step 3. Check the condition (FS 1) in Definition \ref{dn-FS-ga} for $\{\Theta_{K}\}_{K\in \mathcal{D}}$.

Clearly, for all $K\in \mathcal{D}$, $\Theta_{K}\subseteq\operatorname{Id}_{(U, R, \mathcal{F})}$. Conversely, let $F, G\in \mathcal{F}$. If $(F, G)\in \operatorname{Id}_{(U, R, \mathcal{F})}$, then $G\subseteq \overline{R}(F)$. By the proof of Theorem \ref{tm-ToBF-BF}, we have an approximate identity $\{\delta_{K}\}_{K\in \mathcal{D}}$  on $\mathfrak{C}(U, R, \mathcal{F})$,  where  for all $E\in \mathfrak{C}(U, R, \mathcal{F})$, $\delta_{K}(E)=\bigcup\{\overline{R}(M)\mid M\in \mathcal{M}_{K}\mbox{ and } M\subseteq E\}$. Thus $\bigvee_{K\in \mathcal{D}}\delta_{K}(\overline{R}(F))= \overline{R}(F)$. It follows from Lemma \ref{tm-cfga-way}(1) and $R$ being a  preorder that $\overline{R}(F)$ is a compact element of $(\mathfrak{C}(U, R, \mathcal{F}), \subseteq)$. So there exists $J\in \mathcal{D}$ such that $\delta_{J}(\overline{R}(F))= \overline{R}(F)$. Thus $G\subseteq\delta_{J}(\overline{R}(F))$. By Step 2 of Theorem \ref{tm-ToBF-BF}, there exists $M_{F}^{*}\in \mathcal{M}_{J}$ such that
$\delta_{J}(\overline{R}(F))=\overline{R}(M_{F}^{*})\subseteq\overline{R}(F)$.
Hence, $G\subseteq\overline{R}(M_{F}^{*})\subseteq\overline{R}(F)$, showing that $(F, G)\in \Theta_{J}\subseteq\bigcup_{K\in \mathcal{D}}\Theta_{K}$ and $\bigcup_{K\in \mathcal{D}}\Theta_{K}=\operatorname{Id}_{(U, R, \mathcal{F})}$, showing that $\{\Theta_{K}\}_{K\in \mathcal{D}}$ satisfies the condition  (FS 1).

 Step 4. Check that $\{\Theta_{K}\}_{K\in \mathcal{D}}$ satisfies  the condition (FS 2) in Definition \ref{dn-FS-ga}.

  For all $K\in \mathcal{D}$, let  $\mathcal{M}_{K}\subseteq_{fin}\mathcal{F}$ be the one stated in Remark \ref{rk-mk-bfs}(2). Let $F\in \mathcal{F}$.  By the proof of Theorem \ref{tm-ToBF-BF},  there is an $M^{*}\in \mathcal{M}_{K}$ such that $\overline{R}(M^{*})$ is the greatest element in $\{\overline{R}(M)\mid M\in \mathcal{M}_{K}\mbox{ and } M\subseteq \overline{R}(F)\}$. If $ G\in \mathcal{F}$ and $(F, G)\in \Theta_{K}$, then there is $M\in \mathcal{M}_{K}\mbox{ such that } G\subseteq \overline{R}(M)\subseteq \overline{R}(F)$. Noticing that $\overline{R}(M^{*})$ is the greatest element in $\{\overline{R}(M)\mid M\in \mathcal{M}_{K}\mbox{ and } M\subseteq \overline{R}(F)\}$, we have $G\subseteq \overline{R}(M^{*})\subseteq \overline{R}(F)$, showing that $\{\Theta_{K}\}_{K\in \mathcal{D}}$ satisfies the condition (FS 2).

To sum up, $(U, R, \mathcal{F})$ is a topological FS-approxiamtion space.\hfill$\Box$

\begin{tm}\label{tm-BF-ToBF} Let $L$ be a BF-domain, $(K(L), R_{K(L)}, \mathcal{F}_{K(L)})$ the induced CF-approximation space by $L$. Then $(K(L), R_{K(L)}, \mathcal{F}_{K(L)})$ is a topological BF-approxiamtion space.
\end{tm}
\noindent{\bf Proof.} By Remark \ref{GA-Pre-topCF},  $(K(L), R_{K(L)}, \mathcal{F}_{K(L)})$ is a topological CF-approxiamtion space.   By Lemma \ref{lm-BF-eq}, there is  an approximate identity $\{\delta_{i}\}_{i\in I}$  for $L$ consisting of kernel operators with finite range. For all $i\in I$,  we use Im$(\delta_{i})$ to denote the range of $\delta_{i}$.  For any $m\in$ Im$(\delta_{i})$, by Lemma \ref{lm-fis-map} (2), we have that $m=\delta_{i}(m)\ll m$ and $m$ is compact, showing that Im$(\delta_{i})\subseteq K(L)$.

For $H\subseteq_{fin} K(L)$ and $H\ne\emptyset$, it follows from  $\bigvee_{i\in I}\delta_{i}(a)=a$ that there is $i_{a}\in I$ such that $a\leqslant\delta_{i_{a}}(a)$ for all $a\in H\subseteq K(L)$. Clearly, $\delta_{i_{a}}(a)\leqslant a$ and $\delta_{i_{a}}(a)= a$. For $\{\delta_{i_{a}}\mid a\in H\}\subseteq_{fin}\{\delta_{i}\}_{i\in I}$, it follows from the directedness of $\{\delta_{i}\}_{i\in I}$ that there exists $j\in I$ such that $\delta_{i_{a}}\leqslant\delta_{j}$ for all $a\in H$. Thus for all $a\in H$, we have that  $a\geqslant\delta_{j}(a)\geqslant\delta_{i_{a}}(a)=a$, and $a=\delta_{j}(a)$, showing that $H\subseteq\mbox{Im}(\delta_{j})$. Set $\mathcal{M}_{H}=\mathcal{P}(\mbox{Im}(\delta_{j}))\cap\mathcal{F}_{K(L)}$. It follows from $H\subseteq\mbox{Im}(\delta_{j})$ that $\mathcal{P}(H)\cap\mathcal{F}_{K(L)}\subseteq \mathcal{M}_{H}$, showing that $\mathcal{M}_{H}$ satisfies (TB 1).
 To show $\mathcal{M}_{H}$ satisfies (TB 2), let $F\in \mathcal{F}$ and $\mathcal{G}\subseteq\mathcal{M}_{H}$ satisfying $\bigcup\mathcal{G}\subseteq \overline{R}_{K(L)}(F)={\downarrow c_{F}}\cap K(L)$. Take $M=\{\delta_{j}(c_{F})\}\subseteq\mbox{Im}(\delta_{j})$. Clearly, $M\in \mathcal{M}_{H}$. It follows from $\delta_{j}(c_{F})\leqslant c_{F}$ that $\overline{R}_{K(L)}(M)={\downarrow \delta_{j}(c_{F})}\cap K(L)\subseteq{\downarrow c_{F}}\cap K(L)=\overline{R}_{K(L)}(F)$. For all $g\in \bigcup\mathcal{G}$, notice that $\bigcup\mathcal{G}\subseteq\mbox{Im}(\delta_{j})$, we have $\delta_{j}(g)=g$.
It follows from $\bigcup\mathcal{G}\subseteq \overline{R}_{K(L)}(F)={\downarrow c_{F}}\cap K(L)$ that $g\leqslant c_{F}$. Thus $g=\delta_{j}(g)\leqslant\delta_{j}(c_{F})$, showing that $\bigcup\mathcal{G}\subseteq {\downarrow \delta_{j}(c_{F})}\cap K(L)=\overline{R}_{K(L)}(M)$. Thus, we obtain that $\bigcup\mathcal{G}\subseteq\overline{R}_{K(L)}(M)\subseteq\overline{R}_{K(L)}(F)$, showing that $\mathcal{M}_{H}$ satisfies (TB 2).

For $H\subseteq_{fin} K(L)$ and $H=\emptyset$, we have $\mathcal{P}(H)\cap\mathcal{F}_{K(L)}=\emptyset$. Let $\mathcal{M}_{\emptyset}=\mathcal{P}(\mbox{Im}(\delta_{i}))\cap\mathcal{F}_{K(L)}$, for any $i\in I$. Obviously, $\mathcal{M}_{\emptyset}\subseteq_{fin}\mathcal{F}_{K(L)}$ and $\mathcal{M}_{\emptyset}$ satisfies (TB 1). The check of (TB 2) is similar to the case of $H\ne\emptyset$. \hfill$\Box$

\begin{tm}  {\em (Representation Theorem IV: for BF-domains)}   A  poset $(L, \leqslant)$ is  a BF-domain iff there exists a topological BF-approximation space $(U, R, \mathcal{F})$ such that $(\mathfrak{C}(U, R, \mathcal{F}), \subseteq)\cong (L, \leqslant)$.

\end{tm}
\noindent{\bf Proof.} Follows directly from  Proposition \ref{pn-topBF>topFS},  Theorem \ref{rep-BF-topF}, \ref{tm-ToBF-BF} and \ref{tm-BF-ToBF}.\hfill$\Box$

\section{Categorical equivalence between related categories}

 In this section, we first define compositions of suitable CF-approximable relations, and then prove categorical equivalence of related categories.

\begin{dn} Let $(U_{1},  R_{1}, \mathcal{F}_{1})$, $(U_{2},  R_{2}, \mathcal{F}_{2})$ and $(U_{3},  R_{3}, \mathcal{F}_{3})$ be CF-approximation spaces, $\mathrel{\Theta}\subseteq\mathcal{F}_{1}\times\mathcal{F}_{2}$ and $\mathrel{\Upsilon}\subseteq\mathcal{F}_{2}\times\mathcal{F}_{3}$ be CF-approximable relations. Then the composition  $\mathrel{\Upsilon}\circ\mathrel{\Theta}\subseteq\mathcal{F}_{1}\times\mathcal{F}_{3}$ of $\mathrel{\Upsilon}$ and $\mathrel{\Theta}$ is define by that
for any $F_{1}\in \mathcal{F}_{1}, F_{3}\in \mathcal{F}_{3}$, $(F_{1}, F_{3})\in \mathrel{\Upsilon}\circ\mathrel{\Theta}$ iff there exists $F_{2}\in \mathcal{F}_{2}$ satisfying $(F_{1}, F_{2})\in \mathrel{\Theta}$  and $(F_{2}, F_{3})\in \mathrel{\Upsilon}$.
\end{dn}
\begin{pn}\label{pn-comp-cf-app-rel}
Let $(U_{1},  R_{1}, \mathcal{F}_{1})$, $(U_{2},  R_{2}, \mathcal{F}_{2})$ and $(U_{3},  R_{3}, \mathcal{F}_{3})$ be CF-approximation spaces, $\mathrel{\Theta}\subseteq\mathcal{F}_{1}\times\mathcal{F}_{2}$ and $\mathrel{\Upsilon}\subseteq\mathcal{F}_{2}\times\mathcal{F}_{3}$ be CF-approximable relations. Then the composition  $\mathrel{\Upsilon}\circ\mathrel{\Theta}$ is a CF-approximable relation from CF-approximation space $(U_{1}, R_{1}, \mathcal{F}_{1})$ to  $(U_{3}, R_{3}, \mathcal{F}_{3})$.
\end{pn}
\noindent{\bf Proof.}
Since the composition of CF-approximable relation ``$\circ$" is precisely the composition of binary relation, we have that  $\mathrel{\Upsilon}\circ\mathrel{\Theta}$
satisfies the condition (1) in Definition \ref{dn-CF-app}.

To check  that  $\mathrel{\Upsilon}\circ\mathrel{\Theta}$
satisfies the condition (2) in Definition \ref{dn-CF-app},
let $F, F^{\prime}\in \mathcal{F}_{1}$, $H\in  \mathcal{F}_{3}$.
\begin{align*}
&F\subseteq \overline{R_{1}}(F^{\prime}), (F, H)\in \mathrel{\Upsilon}\circ\mathrel{\Theta}\\
\Rightarrow& \exists
 G\in \mathcal{F}_{2}\mbox{ such that }F\subseteq \overline{R_{1}}(F^{\prime}), F\mathrel{\Theta} G\mbox{ and }G\mathrel{\Upsilon} H\\
\Rightarrow& \exists
 G\in \mathcal{F}_{2}\mbox{ such that } F^{\prime}\mathrel{\Theta} G\mbox{ and }G\mathrel{\Upsilon} H \mbox{ (by $\mathrel{\Theta}$ satisfying (2) in  Definition }\ref{dn-CF-app})\\
\Leftrightarrow& (F^{\prime}, H)\in \mathrel{\Upsilon}\circ\mathrel{\Theta}.
\end{align*}
This shows that $\mathrel{\Upsilon}\circ\mathrel{\Theta}$ satisfies the condition (2) in Definition \ref{dn-CF-app}. Similarly, by $\mathrel{\Upsilon}$ satisfying the condition (3) in Definition \ref{dn-CF-app}, we have that $\mathrel{\Upsilon}\circ\mathrel{\Theta}$ satisfies the condition (3) in Definition \ref{dn-CF-app}.

To check  that  $\mathrel{\Upsilon}\circ\mathrel{\Theta}$
satisfies the condition (4) in Definition \ref{dn-CF-app},
let $F\in \mathcal{F}_{1}$,  $H\in \mathcal{F}_{3}$.
$$\begin{array}{lll}
&&(F, H)\in \mathrel{\Upsilon}\circ\mathrel{\Theta}\\
&\Rightarrow&\exists G\in \mathcal{F}_{2}\mbox{ such that }F\mathrel{\Theta} G, G\mathrel{\Upsilon} H\\
&\Rightarrow&\exists F^{\prime}\in \mathcal{F}_{1}, G\in \mathcal{F}_{2}, H^{\prime}\in \mathcal{F}_{3}\mbox{ such that } F^{\prime}\subseteq \overline{R_{1}}(F), F^{\prime}\mathrel{\Theta} G; G\mathrel{\Upsilon} H^{\prime}, H\subseteq \overline{R_{3}}(H^{\prime})\\
&&(\mbox{by Proposition  }\ref{pn1-CF-apprel})\\
&\Rightarrow&\exists F^{\prime}\in \mathcal{F}_{1}, H^{\prime}\in \mathcal{F}_{3}\mbox{ such that }F^{\prime}\subseteq \overline{R_{1}}(F), H\subseteq \overline{R_{3}}(H^{\prime}) \mbox{ and }(F^{\prime}, H^{\prime})\in \mathrel{\Upsilon}\circ\mathrel{\Theta}.

\end{array}$$

To check  that  $\mathrel{\Upsilon}\circ\mathrel{\Theta}$
satisfies the condition (5) in Definition \ref{dn-CF-app},
let $F\in \mathcal{F}_{1}$, $H_{1}, H_{2}\in \mathcal{F}_{3}$.
$$\begin{array}{lll}
&&(F, H_{1})\in  \mathrel{\Upsilon}\circ\mathrel{\Theta}\mbox{ and }(F, H_{2})\in  \mathrel{\Upsilon}\circ\mathrel{\Theta}\\
&\Rightarrow& \exists G_{1}, G_{2}\in \mathcal{F}_{2}\mbox{ such that } F\mathrel{\Theta} G_{1}, F\mathrel{\Theta} G_{2};  G_{1}\mathrel{\Upsilon} H_{1}, G_{2}\mathrel{\Upsilon} H_{2} \\
&\Rightarrow& \exists G_{3}\in \mathcal{F}_{2}\mbox{ such that }G_{1}\cup G_{2}\subseteq \overline{R_{2}}(G_{3}), F\mathrel{\Theta} G_{3}; G_{1}\mathrel{\Upsilon} H_{1}, G_{2}\mathrel{\Upsilon} H_{2}\\
&&(\mbox{by }\mathrel{\Theta}\mbox{ satisfying (5) in Definition }\ref{dn-CF-app})\\
&\Rightarrow&  G_{3}\mathrel{\Upsilon} H_{1}, G_{3}\mathrel{\Upsilon} H_{2}\mbox{ and }F\mathrel{\Theta} G_{3}~~(\mbox{by }\mathrel{\Upsilon}\mbox{ satisfying (2) in Definition }\ref{dn-CF-app})\\
&\Rightarrow& \exists H_{3}\in \mathcal{F}_{3}\mbox{ such that }H_{1}\cup H_{2}\subseteq\overline{R_{3}}(H_{3}), G_{3}\mathrel{\Upsilon} H_{3}\mbox{ and }F\mathrel{\Theta} G_{3}\\
&&(\mbox{by }\mathrel{\Upsilon}\mbox{ satisfying (5) in Definition }\ref{dn-CF-app})\\
&\Rightarrow&\exists H_{3}\in \mathcal{F}_{3}\mbox{ such that }H_{1}\cup H_{2}\subseteq\overline{R_{3}}(H_{3})\mbox{ and } (F, H_{3})\in  \mathrel{\Upsilon}\circ\mathrel{\Theta}.

\end{array}$$

To sum up,  $\mathrel{\Upsilon}\circ\mathrel{\Theta}$ is a CF-approxiamble relation from $(U_{1},  R_{1}, \mathcal{F}_{1})$ to $(U_{3},  R_{3}, \mathcal{F}_{3})$.
\hfill$\Box$
\begin{pn}\label{pn-ide-comp} Let $\mathrel{\Theta}$ be a CF-approximable relation from CF-approximation space $(U_{1}, R_{1}, \mathcal{F}_{1})$ to  $(U_{2}, R_{2}, \mathcal{F}_{2})$. Then $\mathrel{\Theta}\circ \operatorname{Id}_{(U_{1},  R_{1}, \mathcal{F}_{1})}=\operatorname{Id}_{(U_{2},  R_{2}\mathcal{F}_{2})}\circ\mathrel{\mathrel{\Theta}}=\mathrel{\Theta}$, where the identity CF-approximable relation $\operatorname{Id}_{(U_{1},  R_{1}, \mathcal{F}_{1})}$ is defined in Definition \ref{dn-cf-app-rel-id}. \end{pn}
\noindent{\bf Proof.} For all $F\in \mathcal{F}_{1}, G\in \mathcal{F}_{2}$, we have
\begin{align*}
(F, G)\in \mathrel{\Theta}&\Leftrightarrow \exists G^{\prime}\in \mathcal{F}_{2} ~s.t.~ F\mathrel{\Theta} G^{\prime}, G\subseteq \overline{R_2}(G^{\prime})
(\mbox{by Proposition}~\ref{pn1-CF-apprel})\\
&\Leftrightarrow \exists G^{\prime}\in \mathcal{F}_{2} ~s.t. ~F\mathrel{\Theta} G^{\prime}, (G^{\prime}, G)\in \operatorname{Id}_{(U_{2},  R_{2}, \mathcal{F}_{2})}\\
&\Leftrightarrow(F, G)\in \operatorname{Id}_{(U_{2},  R_{2}, \mathcal{F}_{2})}\circ\mathrel{\Theta}.
\end{align*}
This shows that $\operatorname{Id}_{(U_{2},  R_{2}, \mathcal{F}_{2})}\circ\mathrel{\Theta}=\Theta$. Similarly, by Proposition \ref{pn1-CF-apprel}, we have $\Theta\circ \operatorname{Id}_{(U_{1},  R_{1}, \mathcal{F}_{1})}=\Theta$. \hfill$\Box$\\


 By Proposition \ref{pn-comp-cf-app-rel} and \ref{pn-ide-comp},   FS-approximation spaces (resp., strong FS-approximation spaces, topological FS-approximation spaces, topological BF-approximation spaces) as  objects and CF-approximable relations as morphisms
  form a category, denoted by {\bf FS-APPS} (resp.,  {\bf StFS-APPS}, {\bf TopFS-APPS}, {\bf TopBF-APPS}).

  In this paper, we use {\bf FS-DOM}  (resp., {\bf BF-DOM}) to  denote  the category of  FS-domains (resp., BF-domains) and  Scott continuous maps.

Next, we  investigate more relationships between CF-approximable relations and Scott continuous maps, and establish categorical equivalence between  category {\bf FS-APPS} (resp.,  {\bf TopFS-APPS}) and  category {\bf FS-DOM} (resp., {\bf BF-DOM}).

\begin{lm}{\em (\cite{Wu-Xu})}\label{lm-sco-app}
Let $\mathrel{\Theta}$ be  a     CF-approximable relation from     $(U_{1}, R_{1}, \mathcal{F}_{1})$ to  $(U_{2}, R_{2}, \mathcal{F}_{2})$.

$(1)$ Let $f: \mathfrak{C}(U_{1},  R_{1}, \mathcal{F}_{1})\longrightarrow\mathfrak{C}(U_{2},  R_{2}, \mathcal{F}_{2})$ be a Scott continuous map.
Define $\mathrel{\Theta}_{f}\subseteq\mathcal{F}_{1}\times\mathcal{F}_{2}$
such that
\vskip3pt

\centerline{$\forall F\in  \mathcal{F}_{1}, G\in \mathcal{F}_{2}, F\mathrel{\Theta}_{f} G\Leftrightarrow G\subseteq f(\overline{R_{1}}(F)).$}
\vskip5pt

\noindent Then $\mathrel{\Theta}_{f}$ is a CF-approximable relation from $(U_{1},  R_{1}, \mathcal{F}_{1})$ to $(U_{2},  R_{2}, \mathcal{F}_{2})$.

$(2)$  Let $\mathrel{\Theta}$ be a CF-approximable relation from CF-approximation spaces $(U_{1}, R_{1}, \mathcal{F}_{1})$ to  $(U_{2}, R_{2}, \mathcal{F}_{2})$, $f: \mathfrak{C}(U_{1},  R_{1}, \mathcal{F}_{1})\longrightarrow\mathfrak{C}(U_{2},  R_{2}, \mathcal{F}_{2})$  a Scott continuous map.  Then
$\mathrel{\Theta}_{f_{\mathrel{\Theta}}}=\mathrel{\Theta}$ and $f_{\mathrel{\Theta}_{f}}=f$, where  $f_{\mathrel{\Theta}}$ is the one defined in Lemma \ref{tm-cfap-sco}.
\end{lm}
\noindent{\bf Proof.} It is a routine work to check relevant items by  Definition \ref{dn-CF-app}, Lemma \ref{pn2-cf-clo} and \ref{pn3-cf-cl}. \hfill$\Box$

\begin{pn}\label{pn-app-sco}
Let $L_{1}$ and $L_{2}$ be  continuous domains, $(L_{1}, R_{L_{1}},  \mathcal{F}_{L_{1}})$ and $(L_{2}, R_{L_{2}}, \mathcal{F}_{L_{2}})$ the relative induced CF-approximation spaces.

$(1)$ For any CF-approximable relation $\Omega$ from $(L_{1}, R_{L_{1}}, \mathcal{F}_{L_{1}})$ to $(L_{2}, R_{L_{2}}, \mathcal{F}_{L_{2}})$, define $g_{\Omega}: L_{1}\to L_{2}$ such that for  any $x\in L_{1}$,
\vskip3pt

\centerline{$g_{\Omega}(x)=\sup_{L_{2}}(\bigcup\{\overline{R_{L_{2}}}(G)\mid F\in \mathcal{F}_{L_{1}}, G\in \mathcal{F}_{L_{2}}, F\subseteq \dda x \mbox{ and } (F, G)\in\Omega\}.$}
\vskip5pt

\noindent  Then $g_{\Omega}$ is Scott continuous.

$(2)$ Let $g: L_{1}\longrightarrow L_{2}$ be a Scott continuous map, $\Omega$  a CF-approximable relation from $(L_{1}, R_{L_{1}}, \mathcal{F}_{L_{1}})$ to $(L_{2}, R_{L_{2}}, \mathcal{F}_{L_{2}})$. Then $g_{\Omega_{g}}=g$ and $\Omega_{{g}_{\Omega}}=\Omega$, 
where  $\Omega_{g}$ is defined in Proposition \ref{pn-dom-map}.
\end{pn}
\noindent{\bf Proof.} (1) It follows from Lemma \ref{tm-cfap-sco} and Example \ref{tm-CDO-CFGA}(2) that $g_{\Omega}(x)=\sup_{L_{2}}f_{\Omega}(\dda x)$ and $f_{\Omega}(\dda x)\in \mathfrak{C}(L_{2}, R_{L_{2}}, \mathcal{F}_{L_{2}})=\{\dda y\mid y\in L_{2}\}$, where $f_{\Omega}$  is defined in Lemma \ref{tm-cfap-sco}. By continuity of $L_{2}$, we have $f_{\Omega}(\dda x)=\dda y$ for some $y\in L_2$ and $g_{\Omega}$ is well defined. For any directed set $D\subseteq L_{1}$, we have
$$\begin{array}{lll}
g_{\Omega}(\sup_{L_{1}} D)&=&\sup_{L_{2}} f_{\Omega}(\dda \sup_{L_{1}} D)\\
&=&\sup_{L_{2}} f_{\Omega}(\bigcup_{d\in D}\dda d) ~(\mbox{by the continuity of }~L_{1})\\
&=&\sup_{L_{2}} (\bigcup_{d\in D}f_{\Omega}(\dda d))~(\mbox{by Lemma \ref{pn2-cf-clo}(3), \ref{tm-cfap-sco} and Example \ref{tm-CDO-CFGA}(2) })\\
&=&\sup_{L_{2}} \{\sup_{L_{2}}f_{\Omega}(\dda d)\mid d\in D\}\\
&=&\sup_{L_{2}} \{g_{\Omega}(d) \mid d\in D\}.
\end{array}$$
This shows that $g_{\Omega}$ is Scott continuous.

(2) For every $x\in L_{1}$ and  $\Omega_g$ defined in Proposition \ref{pn-dom-map}, we have
$$\begin{array}{lll}
g_{\Omega_{g}}(x)&=&\sup_{L_{2}} f_{\Omega_g}(\dda x)\\
&=&\sup_{L_{2}}(\bigcup\{\overline{R_{L_{2}}}(G)\mid F\in \mathcal{F}_{L_{1}}, G\in \mathcal{F}_{L_{2}}, F\subseteq \dda x ~\mbox{and}~ (F, G)\in\Omega_g\}\\
&=&\sup_{L_{2}}(\bigcup\{\dda t\mid s\in L_{1}, t\in L_{2}, s\ll x ~\mbox{and}~ t\ll g(s)\})\\
&=&\sup_{L_{2}} \{ t\mid s\in L_{1}, t\in L_{2}, s\ll x ~\mbox{and}~ t\ll g(s)\}\\
&=&\sup_{L_{2}} \{ g(s)\mid s\in L_{1}~\mbox{and}~s\ll x \}~(\mbox{by the continuity of}~L_{2}~)\\
&=&g(x).~(\mbox{by the Scott continuity of}~g)
\end{array}$$
This shows that $g_{\Omega_{g}}=g$.

 For any $F\in \mathcal{F}_{L_{1}}$  and $G\in \mathcal{F}_{L_{2}}$, then
$$\begin{array}{lll}
(F, G)\in \Omega_{{g}_{\Omega}}&\Longleftrightarrow& c_{G}\ll g_{\Omega}(c_{F})=\sup_{L_2} f_{\Omega}(\dda c_{F})\\
&\Longleftrightarrow& G\subseteq f_{\Omega}(\dda c_{F})~ (\mbox{by}~f_{\Omega}(\dda x)\in \mathfrak{C}(L_{2}, R_{L_{2}}, \mathcal{F}_{L_{2}})=\{\dda y\mid y\in L_{2}\})\\
&\Longleftrightarrow& G\subseteq \bigcup\{\overline{R_{L_{2}}}(K)\mid H\in \mathcal{F}_{L_{1}}, K\in \mathcal{F}_{L_{2}}, H\subseteq \dda c_F ~\mbox{and}~ (H, K)\in\Omega\}\\
&\Longleftrightarrow& (\exists H\in \mathcal{F}_{L_{1}},  K\in \mathcal{F}_{L_{2}})(H\subseteq \dda c_{F}=\overline{R_{L_{1}}}(F), G\subseteq \overline{R_{L_{2}}}(K) ~\mbox{and}~ (H, K)\in \Omega)\\
&\qquad& \ (\mbox{by Proposition }\ref{pn-ovRF-direct} )\\
&\Longleftrightarrow& (F, G)\in \Omega.~(\mbox{by Proposition \ref{pn1-CF-apprel}})
\end{array}$$
 This shows that $\Omega_{{g}_{\Omega}}=\Omega$.
\hfill$\Box$

\begin{tm} \label{lm-FSGA-FSDOM} Let $(U, R, \mathcal{F})$ be an FS-approximation space, $(V, Q,  \mathcal{G})$ the FS-approximation space induced by the FS-domain $(\mathfrak{C}(U, R, \mathcal{F}), \subseteq)$. Then $(U, R, \mathcal{F})\cong(V, Q,  \mathcal{G})$ in  {\bf FS-APPS}.
\end{tm}
\noindent{\bf Proof.}
By Example \ref{tm-CDO-CFGA}, we have  that $V= \mathfrak{C}(U, R, \mathcal{F})$, $\mathcal{G}=\{\mathcal{C}\subseteq\mathfrak{C}(U, R, \mathcal{F})\mid \mathcal{C} ~has~a~top~element\}$, and $Q$ is the way-below relation  $\ll$ in FS-domain $(\mathfrak{C}(U, R, \mathcal{F}), \subseteq)$.
Define $\mathrel{\Upsilon}\subseteq\mathcal{F}\times\mathcal{G}$ and $\Omega\subseteq\mathcal{G}\times\mathcal{F}$ such that
 \vskip3pt

 \centerline{$\forall F\in  \mathcal{F},  \mathcal{C}\in \mathcal{G}, (F, \mathcal{C})\in \mathrel{\Upsilon}\Leftrightarrow E_{\mathcal{C}}\ll \overline{R}(F), \indent (\mathcal{C}, F)\in \Omega\Leftrightarrow F\subseteq E_{\mathcal{C}},$}
 \vskip5pt
\noindent where  $E_{\mathcal{C}}$ denote the top element in $\mathcal{C}$.

 We prove firstly that $\mathrel{\Upsilon}$ is  a CF-aproximable relation from $(U, R, \mathcal{F})$ to $(V, Q,  \mathcal{G})$. It is a routine work to check that $\mathrel{\Upsilon}$ satisfies the conditions (1), (2), (3) and (5) in Definition \ref{dn-CF-app}. To check that  $\mathrel{\Upsilon}$ satisfies (4) in Definition \ref{dn-CF-app}, let $ F\in  \mathcal{F},  \mathcal{C}\in \mathcal{G}$. Then we have
\begin{align*}
(F, \mathcal{C})\in \mathrel{\Upsilon}  &\Longleftrightarrow E_{\mathcal{C}}\ll\overline{R}(F) \\
&\Longleftrightarrow \exists E_{1}, E_{2}\in \mathfrak{C}(U, R, \mathcal{F})\mbox{ such that } E_{\mathcal{C}}\ll E_{1}\ll E_{2}\ll \overline{R}(F)~(\mbox{by Lemma}~\ref{pn-int})\\
&\Longleftrightarrow \exists E_{1}, E_{2}\in \mathfrak{C}(U, R, \mathcal{F}), F^{\prime}\in \mathcal{F} \mbox{ such that } E_{\mathcal{C}}\ll E_{1}\ll E_{2}\subseteq \overline{R}(F^{\prime}), F^{\prime}\subseteq\overline{R}(F)\\
&~\qquad (\mbox{by Lemma}~\ref{tm-cfga-way})\\
&\Longrightarrow  \exists E_{1}\in \mathfrak{C}(U, R, \mathcal{F}), F^{\prime}\in \mathcal{F}\mbox{ such that } E_{\mathcal{C}}\ll E_{1}, E_{1}\ll \overline{R}(F^{\prime}), F^{\prime}\subseteq\overline{R}(F)\\
&\Longrightarrow  \exists F^{\prime}\in \mathcal{F}, \{E_{1}\}\in \mathcal{G}\mbox{ such that } F^{\prime}\subseteq\overline{R}(F), (F^{\prime}, \{E_{1}\})\in\mathrel{\Upsilon},  \mathcal{C}\subseteq \overline{Q}(\{ E_{1}\})=\dda \{E_1\}.
\end{align*}
This shows that that  $\mathrel{\Upsilon}$ satisfies (4) in Definition \ref{dn-CF-app}.

 We then prove  that $\Omega$ is  a CF-aproximable relation from $(V, Q,  \mathcal{G})$ to $(U, R, \mathcal{F})$. It is a routine work to check that $\Omega$ satisfies the conditions (1), (2), (3) and (5) in Definition \ref{dn-CF-app}. To check that  $\Omega$ satisfies (4) in Definition \ref{dn-CF-app}, let $\mathcal{C}\in \mathcal{G},  F\in  \mathcal{F}$. Then we have
\begin{align*}
(\mathcal{C}, F)\in \Omega&\Leftrightarrow F\subseteq E_{\mathcal{C}}\in \mathfrak{C}(U, R, \mathcal{F})\\
&\Rightarrow  \exists F_{1}\in \mathcal{F}\mbox{ such that } F\subseteq \overline{R}( F_{1}),  F_{1}\subseteq E_{\mathcal{C}}~(\mbox{by Definition}~\ref{dn-cf-cl})\\
&\Rightarrow  \exists F_{1}, F_{2}\in \mathcal{F}\mbox{ such that } F\subseteq \overline{R}( F_{1}),  F_{1}\subseteq\overline{R}( F_{2}) \mbox{ and }  F_{2}\subseteq E_{\mathcal{C}}~(\mbox{by Definition}~\ref{dn-cf-cl})\\
&\Rightarrow\exists F_{1}\in \mathcal{F},  \{\overline{R}(F_{2})\}\in  \mathcal{G}\mbox{ such that } F\subseteq \overline{R}( F_{1}), (\{\overline{R}( F_{2})\}, F_{1})\in \Omega\mbox{ and }\overline{R}(F_{2})\ll E_{\mathcal{C}}\\
&~\quad(\mbox{by Lemma}~\ref{tm-cfga-way}\mbox{(2)})\\
&\Rightarrow \exists \{\overline{R}(F_{2})\}\in  \mathcal{G}, F_{1}\in \mathcal{F}\mbox{ such that }\{\overline{R}(F_{2})\}\subseteq \overline{Q} (\mathcal{C}), (\{\overline{R}( F_{2})\}, F_{1})\in \Omega,\ F\subseteq \overline{R}( F_{1}).
\end{align*}
This shows that   $\Omega$ satisfies (4) in Definition \ref{dn-CF-app}.

Next we prove  $\mathrel{\Omega}\circ\mathrel{\Upsilon}=\operatorname{Id}_{(U, R, \mathcal{F})}$ and $ \mathrel{\Upsilon}\circ\Omega=\operatorname{Id}_{(V, Q, \mathcal{G})}$.
Let $F, G\in \mathcal{F}$, $ \mathcal{C}_{1}, \mathcal{C}_{2}\in \mathcal{G}$. Then
\begin{align*}
(F, G)\in \operatorname{Id}_{(U, R, \mathcal{F})}&\Leftrightarrow G\subseteq \overline{R}(F)\\
&\Leftrightarrow \exists F_{1}\in \mathcal{F},G\subseteq\overline{R}(F_{1}),  F_{1}\subseteq\overline{R}(F)~(\mbox{by Definition}~\ref{dn-cf-ga})\\
&\Leftrightarrow \exists \mathcal{C}=\{\overline{R}(F_{1})\}\in \mathcal{G}, G\subseteq E_{ \mathcal{C}},  E_{ \mathcal{C}}\ll\overline{R}(F)~(\mbox{by Lemma}~\ref{tm-cfga-way}(2))\\
&\Leftrightarrow \exists \mathcal{C}\in \mathcal{G}, (F, \mathcal{C})\in \mathrel{\Upsilon},  (\mathcal{C}, G)\in \Omega \mbox { (by definitions of }\Omega\mbox { and }\mathrel{\Upsilon}) \\
&\Leftrightarrow (F, G)\in \mathrel{\Omega}\circ\mathrel{\Upsilon},
\end{align*}
\begin{align*}
(\mathcal{C}_{1}, \mathcal{C}_{2})\in \operatorname{Id}_{(V, Q, \mathcal{G})}&\Leftrightarrow \mathcal{C}_{2}\subseteq \overline{Q}(\mathcal{C}_{1})\\
&\Leftrightarrow E_{\mathcal{C}_{2}}\ll E_{\mathcal{C}_{1}}\\
&\Leftrightarrow \exists E\in \mathfrak{C}(U, R, \mathcal{F}), E_{\mathcal{C}_{2}}\ll E\ll E_{\mathcal{C}_{1}}~(\mbox{by Lemma}~\ref{pn-int})\\
&\Leftrightarrow \exists F\in  \mathcal{F}, F\subseteq E_{\mathcal{C}_{1}},  E_{\mathcal{C}_{2}}\ll \overline{R}(F)~(\mbox{by  taking}~E=\overline{R}(F)\mbox{ and Lemma}~\ref{tm-cfga-way})\\
&\Leftrightarrow \exists F\in  \mathcal{F}, (\mathcal{C}_{1}, F)\in\Omega , (F, \mathcal{C}_{2})\in \mathrel{\Upsilon} \mbox { (by definitions of }\Omega \mbox { and } \mathrel{\Upsilon}) \\
&\Leftrightarrow  (\mathcal{C}_{1},   \mathcal{C}_{2})\in \mathrel{\mathrel{\Upsilon}}\circ\mathrel{\Omega}.
\end{align*}
These show  that $(U, R, \mathcal{F})\cong(V, Q,  \mathcal{G})$ in  {\bf FS-APPS}.\hfill$\Box$

\begin{lm}{\em(\cite{Well})}\label{lm-cat-eqv} Let $\mathcal{C}, \mathcal{D}$ be two categories. If there is a functor $\Phi: \mathcal{C}\longrightarrow \mathcal{D}$ such that

$(1)$ $\Phi$ is  full, namely, for all $A, B\in ob(\mathcal{C})$, $g\in Mor_{\mathcal{D}}(\Phi(A), \Phi(B))$,  there is  $f\in Mor_{\mathcal{C}}(A, B)$ such that $\Phi(f)=g$;

$(2)$ $\Phi$ is faithful, namely, for all $A, B\in ob(\mathcal{C})$, $f, g\in Mor_{\mathcal{C}}(A, B)$, if $f\neq g$, then $\Phi(f)\neq \Phi(g)$;

$(3)$ for all $B\in ob(\mathcal{D})$, there is $A\in ob(\mathcal{C})$ such that $\Phi(A)\cong B$,

\noindent then $\mathcal{C}$ and $\mathcal{D}$ are  equivalent.
\end{lm}

\begin{tm}\label{cat-FSGA-FSDO}
The categories {\bf FS-APPS}  and  {\bf FS-DOM} are equivalent.
\end{tm}
\noindent{\bf Proof.}
  Define  $\Phi:  \mbox{{\bf FS-DOM}}\longrightarrow ${\bf FS-APPS}, such that
\vskip3pt

\centerline{$\forall L\in ob(\mbox{{\bf FS-DOM}}), \Phi(L)=(L, R_{L}, \mathcal{F}_{L});  \quad \forall g\in Mor(\mbox{{\bf FS-DOM}}), \Phi(g)=\Omega_{g},$}
\vskip5pt
\noindent where the $\Omega_{g}$ is defined in Proposition \ref{pn-dom-map}. By  Proposition \ref{pn-dom-map} and  Theorem \ref{tm-FSDO-FSGA}, we know $\Phi$ is well-defined.  We then  prove that $\Phi$ is a functor from \mbox{{\bf FS-DOM}} to {\bf FS-APPS}.
Let $L$ be an FS-domain, $id_{L}$ the identity on $L$. Then for all $F, G\in \mathcal{F}_{L}$, we have
$$\begin{array}{lll}
(F, G)\in  \Phi(id_{L})&\Longleftrightarrow&c_{G}\ll id_L(c_{F})=c_F\\
&\Longleftrightarrow& G\subseteq \overline{R_{L}}(F)=\dda c_F\\
&\Longleftrightarrow& (F, G)\in\operatorname{Id}_{(L, R_{L}, \mathcal{F}_{L})},
\end{array}$$
\noindent where $c_F$ is the top element of $F$. This shows that  $\Phi(id_{L})=\operatorname{Id}_{(L, R_{L}, \mathcal{F}_{L})}$.

Let $L_{1}$,  $L_{2}$, $L_{3}$ be  FS-domains,
 $h: L_{1}\longrightarrow L_{2}$ and $g: L_{2}\longrightarrow L_{3}$ Scott continuous maps. Then for all $F\in \mathcal{F}_{L_{1}}$, $G\in \mathcal{F}_{L_{3}}$, we have
$$\begin{array}{lll}
(F, G)\in  \Phi(g)\circ\Phi(h)&\Longleftrightarrow&\exists K\in \mathcal{F}_{L_{2}}\mbox{ such that } (F, K)\in \Phi(h)\mbox{ and }(K ,G)\in \Phi(g)\\
&\Longleftrightarrow& \exists K\in \mathcal{F}_{L_{2}}, c_{K}\ll h(c_{F}), c_{G}\ll g(c_{K}) \\
&\Longleftrightarrow&  c_{G}\ll  g(h(c_{F})) (``\Leftarrow"\mbox{ follows from } g(h(c_{F}))=\sup\{g(x)\mid x\ll h(c_{F})\}\\
&\Longleftrightarrow& (F, G)\in  \Phi(g\circ h).
\end{array}$$
This shows that $\Phi(g)\circ\Phi(h)= \Phi(g\circ h)$. Thus $\Phi$ is a functor from \mbox{{\bf FS-DOM}} to {\bf FS-APPS}.

It follows from Proposition \ref{pn-app-sco} that $\Phi$ is full and faithful. By Theorem \ref{lm-FSGA-FSDOM}, we see  that $\Phi$ satisfies the condition (3) in Lemma \ref{lm-cat-eqv}. Thus, \mbox{{\bf FS-DOM}} and {\bf FS-APPS} are equivalent. \hfill$\Box$

\begin{tm}
The categories  {\bf TopFS-APPS}  and  {\bf BF-DOM} are equivalent.
\end{tm}
\noindent{\bf Proof.}
Define $\Psi:$ {\bf TopFS-APPS}$\longrightarrow$ {\bf BF-DOM} such that
\vskip3pt

\centerline{$\forall (U, R, \mathcal{F})\in ob(\mbox{{\bf FS-DOM}}),\ \ \Psi((U, R, \mathcal{F}))=(\mathfrak{C}(U, R, \mathcal{F}), \subseteq);$}
\vskip3pt

 \centerline{$\forall \Theta\in Mor(\mbox{{\bf TopFS-APPS}}),\   \Psi(\Theta)=f_{\Theta},$}
\vskip5pt
\noindent where $f_{\Theta}$ is defined in  Lemma \ref{tm-cfap-sco}. By Lemma \ref{tm-cfap-sco} and Theorem \ref{rep-BF-topF},  we know that $\Psi$ is well-defined. Similar to the proof of Theorem \ref{cat-FSGA-FSDO}, it is routine to check that $\Psi$ is a functor. It follows from Theorem \ref{rep-BF-topF} and Lemma \ref{lm-sco-app} that $\Psi$ satisfies  all the conditions in Lemma \ref{lm-cat-eqv}, showing that categories  {\bf TopFS-APPS}  and  {\bf BF-DOM} are equivalent.
\hfill$\Box$\\

 Similarly, the category {\bf StFS-APPS} of strong FS-approximation spaces and CF-approximable relations is equivalent to category {\bf FS-DOM},  and the category {\bf TopBF-APPS} of topological BF-approximation spaces and CF-approximable relations is equivalent to category {\bf BF-DOM}. We leave the details of these to the interested reader.\\

\noindent {\bf Declaration of competing interest}

The authors declare that they have no known competing financial interests or personal relationships that could have appeared to influence the work reported in this paper.


\end{document}